\documentclass[11pt,reqno]{amsart}

\usepackage{amssymb}
\usepackage{amsthm}
\usepackage{amsmath,amsfonts,amssymb}
\usepackage{graphics,color}
\usepackage{enumerate}

\newtheorem{theorem}{Theorem}[section]
\newtheorem{lemma}[theorem]{Lemma}
\newtheorem{proposition}[theorem]{Proposition}

\newtheorem{definition}[theorem]{Definition\rm}
\newtheorem{remark}{Remark}

\newcommand*{\R}{\ensuremath{\mathbb{R}}}
\renewcommand*{\S}{\ensuremath{\mathcal{S}}}
\newcommand*{\N}{\ensuremath{\mathbb{N}}}
\newcommand*{\Z}{\ensuremath{\mathbb{Z}}}

\newcommand*{\supp}{\ensuremath{\mathrm{supp\,}}}

\renewcommand*{\div}{\ensuremath{\mathrm{div\,}}}
\newcommand*{\tr}{\ensuremath{\mathrm{tr\,}}}

\newcommand{\eps}{\varepsilon}
\def\rn#1{\mathbb{R}^{#1}}

\newcommand*{\trs}{\;{\scriptstyle{\bigcirc}}\;}

\begin{document}

\title[On weak solutions of the Euler equations]
{On admissibility criteria for weak solutions 
of the Euler equations}

\author{Camillo De Lellis}
\address{Institut f\"ur Mathematik, Universit\"at Z\"urich, CH-8057 Z\"urich}
\email{camillo.delellis@math.unizh.ch}

\author{L\'aszl\'o Sz\'ekelyhidi Jr.}
\address{Hausdorff Center for Mathematics, Universit\"at Bonn, D-53115 Bonn}
\email{laszlo.szekelyhidi@hcm.uni-bonn.de}

\begin{abstract}
We consider solutions to the Cauchy problem
for the incompressible Euler equations satisfying
several additional requirements, like the global
and local energy inequalities. Using some techniques
introduced in an earlier paper we show that, for 
some bounded compactly supported initial data, none
of these admissibility criteria singles out a unique
weak solution.

As a byproduct we show bounded initial data
for which admissible solutions 
to the $p$--system of isentropic gas dynamics in Eulerian
coordinates are not unique in more than one space dimension.
\end{abstract}

\maketitle

%{\bf Plan of the paper:}

%\begin{itemize}

%\item {\bf Introduction: two main theorems}
%\item {\bf An overview on the different notions of admissible solutions}
%\begin{itemize}
%\item {\bf Weak and strong energy inequalities}
%\item {\bf Local energy inequality}
%\item {\bf Dissipative solutions and uniqueness of smooth solutions}
%\item {\bf Conservation laws}
%\end{itemize}
%\item {\bf Statement of the results}
%\begin{itemize}
%\item {\bf Nonuniqueness for Euler}
%\item {\bf A criterion for the existence of wild solutions}
%\item {\bf Nonuniqueness for isentropic gas--dynamics} 
%\end{itemize}  
%\item {\bf Proof: Part I, the criterion}
%\item {\bf Proof: Part II}
%\item {\bf The system of isentropic gas dynamics}
%\item {\bf Appendices}
%\end{itemize}

\section{Introduction}

In this paper we consider the Cauchy problem for the incompressible
Euler equations in $n$ space dimensions, $n\geq 2$,
\begin{equation}\label{euler}
\left\{
\begin{array}{l}
\partial_tv+\div(v\otimes v)+\nabla p \;=\; 0,\\
\div v\;=\;0,\\
v(x,0)\;=\;v^0(x),  
\end{array}\right.
\end{equation}
where the initial data $v^0$ satisfies 
the compatibility condition
\begin{equation}\label{initial}
\div v^0\;=\; 0\, .
\end{equation}
A vector field $v\in L^2_{loc}(\R^n\times ]0,\infty[)$ is a 
{\em weak solution} of 
\eqref{euler} if $v(\cdot,t)$ is weakly divergence-free for 
almost every $t>0$, and
\begin{equation}\label{distrib}
\int_0^{\infty}\int_{\R^n}
[v \cdot \partial_t \varphi + \langle v\otimes v, 
\nabla \varphi\rangle] \, dx\, dt \;+\; 
\int_{\R^n} v^0(x)\varphi(x, 0)\, dx=0
\end{equation}
for every test function 
$\varphi\in C_c^{\infty}(\R^n\times [0,\infty[; \R^n)$ 
with ${\rm div}\, \varphi = 0$.
It is well--known that then the pressure is determined up to
a function depending only on time (see \cite{TemamBook}).

In his pioneering work \cite{Scheffer93}
V.~Scheffer showed that weak
solutions to the 2--dimensional Euler equations
are not unique. In particular Scheffer constructed a
nontrivial weak solution which is compactly supported in space and time,
thus disproving uniqueness for \eqref{euler} even when $v^0=0$. 
A simpler construction was later proposed by A.~Shnirelman in
\cite{Shnirelman1}.

In a recent paper \cite{us}, we have shown how the general framework
of convex integration \cite{DacorognaMarcellini97,MS99,Bernd,SychevFew} 
combined with Tartar's programme on oscillation 
phenomena in conservation laws \cite{Tartar79} (see also \cite{KMS02} for an overview) can be applied to \eqref{euler}. In this way, one can easily recover Scheffer's and 
Shnirelman's counterexamples in all dimensions
and with bounded velocity and pressure. Moreover, the construction
yields as a simple corollary the existence of
energy--decreasing solutions, thus recovering another
groundbreaking result of Shnirelman
\cite{Shnirelmandecrease}, again with the additional features that
our examples have bounded velocity and pressures and
can be shown to exist in any dimension. 

The results so far left open the question of whether
one might achieve the uniqueness of weak solutions by imposing 
a form of the energy inequality. Our primary purpose in this note
is to address this issue. More precisely we prove the following theorem
(for the relevant definitions of weak, strong and local 
energy inequalities, we refer to Sections 
\ref{s:weakstrongenergy} and \ref{s:localenergy}).

\begin{theorem}\label{t:main}
Let $n\geq 2$. There exist bounded and compactly supported divergence--free vector fields $v^0$
for which there are 
\begin{itemize}
\item[(a)] infinitely many weak solutions of \eqref{euler}
satisfying both the strong and the local energy {\em equalities};
\item[(b)] weak solutions of \eqref{euler}
satisfying the strong energy inequality but not the 
energy equality;
\item[(c)] weak solutions of
\eqref{euler} satisfying the weak energy inequality but not
the strong energy inequality.
\end{itemize}
\end{theorem}

Our examples display very wild behavior, 
such as dissipation of the energy and high--frequency oscillations. We will 
refer to them as {\em wild solutions}. 
A natural question is to characterize the set of initial data $v^0$ to which such wild solutions exist, i.e.~the set of initial data for which Theorem \ref{t:main} holds. The core of this note is devoted to a first characterization in Proposition \ref{p:main} of such "wild" initial data, in terms of the existence of a suitable subsolution. An important point is that - in contrast to the constructions in \cite{Scheffer93,Shnirelman1,us} - for weak solutions satisfying the energy inequality there are nontrivial constraints on $v^0$. For example $v^0$ cannot be smooth (see Section \ref{s:dissipative}). We give a direct construction of wild initial data in Section 5, but for example we were unable to decide the following question:\footnote{as usual $H(\R^n)$ denotes the set of solenoidal vector fields in $L^2(\R^n)$.}
\begin{equation*}
\textrm{is the set of wild initial data dense in $H(\R^n)$?}
\end{equation*}
A related question is to estimate the maximal dissipation rate possible for wild solutions for a given initial data.

As a byproduct of our analysis we prove a similar non--uniqueness
result for the $p$--system of isentropic gas dynamics in Eulerian 
coordinates, the oldest
hyperbolic system of conservation laws.
The unknowns of the system, which consists of $n+1$ equations,
are the density $\rho$ and the velocity $v$ of the gas:
\begin{equation}\label{e:psistema}
\left\{\begin{array}{l}
\partial_t \rho + {\rm div}_x (\rho v) \;=\; 0\\
\partial_t (\rho v) + {\rm div}_x (\rho v\otimes v) + \nabla [ p(\rho)]\;=\; 0\\
\rho (0, \cdot)\;=\; \rho^0\\
v (0, \cdot)\;=\; v^0\, 
\end{array}\right.
\end{equation}
(cf.~(3.3.17) in \cite{DafermosBook} and Section 1.1 of \cite{SerreBook} p7).
The pressure $p$ is a function of $\rho$, which is determined from
the constitutive thermodynamic relations of the gas in question
and satisfies the assumption $p' >0$. A typical example is $p (\rho )= k \rho^\gamma$,
with constants $k>0$ and $\gamma>1$,
which gives the constitutive relation for a polytropic gas
(cf.~(3.3.19) and (3.3.20) of \cite{DafermosBook}).
Weak solutions of \eqref{e:psistema} are bounded functions in $\R^n$, which solve
it in the sense of distributions. Admissible solutions have to satisfy an additional
inequality, coming from the conservation law for the energy of the system.
For the precise definition we refer to Section 2.4.

\begin{theorem}\label{t:psistema}
Let $n\geq 2$. Then, for any given function $p$, 
there exist bounded initial data $(\rho^0, v^0)$ with
$\rho^0\geq c>0$ for which there are infinitely many bounded admissible
solutions $(\rho, v)$ of \eqref{e:psistema} with $\rho\geq c>0$.
\end{theorem}

\begin{remark} In fact, all the solutions constructed in our proof of Theorem
\ref{t:psistema} satisfy the energy {\em equality}. They are therefore also
entropy solutions of the full compressible Euler system (see for instance example
(d) of Section 3.3 of \cite{DafermosBook}) and they show nonuniqueness in this case
as well. This failure of uniqueness was suggsted by Elling in \cite{elling},
although the arguments leading him to this suggestion are completely unrelated
to our setting.
\end{remark}

In fact the same result also holds for the full compressible Euler system, since the solutions we construct satisfy the energy equality, hence there is no entropy production at all.

The paper is organized as follows. Section 2 contains a survey of
several admissibility conditions for \eqref{euler} and the definition
of admissible solutions for \eqref{e:psistema}. Section 3 states
a general criterion on the existence of wild solutions to
\eqref{euler} for a given initial data, in Proposition
\ref{p:main}. 

In Section 4, forming the central part of the paper, we prove Proposition \ref{p:main} by developing a variant of the "Baire category method"
for differential inclusions which is applicable to evolution equations
in the space $C\bigl([0,\infty[;L^2_w(\R^n)\bigr)$ (see below). The Baire
category method has been developed in \cite{BressanFlores,DacorognaMarcellini97,Bernd,SychevFew}, 
and in \cite{us} we applied it to \eqref{euler}. These techniques do not yield solutions
which are weakly continuous in time - a property that is needed in connection with the strong form of the energy inequality. Of course the constructive method is easy to modify to yield
such solutions, but Baire category techniques have the advantage of showing very clearly the arbitrariness
in each step of the construction, by exhibiting infinitely many solutions at the same time. The main point is to find a functional setup in which the points of continuity of a Baire-1 map coincides with solutions of the differential inclusion in the space $C\bigl([0,\infty[;L^2_w(\R^n)\bigr)$. 

In Section 5 we construct initial data meeting the requirements
of Proposition \ref{p:main}, see Proposition \ref{p:initialdata}.
Finally, in Section 6 we prove the non--uniqueness theorems
\ref{t:main} and \ref{t:psistema} using Proposition \ref{p:main}
and Proposition \ref{p:initialdata}.

\section{An overview of the different notions of admissibility}

In this section we discuss various admissibility 
criteria for weak solutions
which have been proposed in the literature. 

%\subsection{Boundary conditions}
%In this paper we will consider also the Euler equations
%in an open domain $\Omega\subset \R^n$ when its boundary is 
%sufficiently regular ($C^1$ suffices). In this case we will
%always assume that $v\cdot \nu |_{\partial \Omega} = 0$ 
%in the weak sense. This simply means that if we extend
%$v$ to be $0$ outside $\Omega$, the extension is a weak solution
%of Euler on the whole $\R^n$ (in the sense that \eqref{distrib} holds). 

%Below we will state all definitions and theorems on further admissibility
%criteria for solutions on the whole space: in the case of 
%solutions on domains $\Omega$ one gets equivalent definitions by simply
%extending the solutions trivially outside $\Omega$.

\subsection{Weak and strong energy inequalities}\label{s:weakstrongenergy}
All the admissibility criteria considered so
far in the literature are motivated by approximating \eqref{euler}
with the Navier Stokes equations. We therefore consider
the following vanishing viscosity approximation of \eqref{euler}
\begin{equation}\label{navierstokes}
\left\{
\begin{array}{l}
\partial_tv+\div(v\otimes v)+\nabla p \;=\; \nu \Delta v\\
\div v\;=\;0\\
v(x,0)\;=\;v^0(x)\, ,
\end{array}\right.
\end{equation}
where the parameter $\nu$ is positive but small.
The \emph{weak formulation} of \eqref{navierstokes}, which makes sense
for any $v\in L^2_{loc}(\R^n\times]0,\infty[)$, is the following:
$v(\cdot,t)$ is weakly divergence-free for almost every $t>0$, and
\begin{equation}\label{distribNS}
\int_0^{\infty}\int_{\R^n} 
\Bigl[ v \cdot (\partial_t \varphi + \nu\Delta\varphi)+ \langle v\otimes v, 
\nabla \varphi\rangle\Bigr] \,dx\,dt\;+\; 
\int_{\R^n} v^0(x)\varphi(x, 0)\, dx=0
\end{equation}
for every test function 
$\varphi\in C_c^{\infty}(\R^n\times [0,\infty[; \R^n)$ 
with ${\rm div}\, \varphi = 0$.

\medskip

For smooth solutions, we can multiply \eqref{euler} 
and \eqref{navierstokes}
by $v$ and derive corresponding partial differential equations 
for $|v|^2$, namely
\begin{equation}\label{LocEnEul}
\partial_t \frac{|v|^2}{2} + \div \left(v \left(\frac{|v|^2}{2} + p\right)\right)
\;=\; 0
\end{equation}
and
\begin{equation}\label{LocEnNS}
\partial_t \frac{|v|^2}{2} + \div \left(v \left(\frac{|v|^2}{2} + p\right)\right)
\;=\; \nu \Delta \frac{|v|^2}{2} - \nu |\nabla v|^2\, .
\end{equation}
Recall that \eqref{euler}
and \eqref{navierstokes} model the movements of ideal 
incompressible fluids.
If we assume that the constant density of the 
fluid is normalized
to $1$, then $|v|^2/2$ is the energy density
and \eqref{LocEnEul} and \eqref{LocEnNS} 
are simply the laws of conservation of the 
energy, in local form. 
%Therefore we will call these equalities 
%{\em local energy equalities}. 

Integrating \eqref{LocEnEul} and \eqref{LocEnNS} in 
time and space and assuming that $p$ and $v$ are decaying sufficiently
fast at infinity, we deduce formally the following identities:
\begin{equation}\label{EnEul}
\frac{1}{2} \int_{\R^n} |v|^2 (x,t)\, dx 
\;=\; \frac{1}{2} \int_{\R^n} |v|^2 (x,s)\, ds \qquad \textrm{for all $s,t\geq 0$,}
\end{equation}
\begin{equation}\label{EnNS}
\frac{1}{2}\int_{\R^n} |v|^2 (x,t)\, dx 
\;=\; \frac{1}{2}\int_{\R^n} |v|^2 (x,s)\, dx - \nu \int_s^t \int_{\R^n}
|\nabla v|^2 (x,\tau)\, dx\, d\tau\, .
\end{equation}

The celebrated result of J.~Leray \cite{leray}
(see \cite{Galdi} for a modern introduction)
shows the existence of weak solutions
to \eqref{navierstokes} which satisfy a relaxed 
version of \eqref{EnNS}.

\begin{theorem}[Leray]\label{t:leray}
Let $v^0\in L^2(\R^n)$ be a divergence--free vector field.
Then there exists $v\in L^{\infty}([0,\infty[;L^2(\R^n))$ with $\nabla v\in L^2(\R^n\times]0,\infty[)$
such that $v(\cdot,t)$ is weakly divergence-free and \eqref{distribNS} holds  for all $t>0$. Moreover,
\begin{equation}\label{e:leray1}
\begin{split}
\frac{1}{2}\int_{\R^n} |v|^2 (x,t)\, dx 
\;\leq\; \frac{1}{2}\int_{\R^n} |v^0|^2 (x)\, dx- \nu \int_0^t \int_{\R^n}
&|\nabla v|^2 (x,\tau)\, dx\, d\tau\\
&\textrm{ for every $t>0$,}
\end{split}
\end{equation}
and more generally
\begin{equation}\label{e:leray2}
\begin{split}
\frac{1}{2}\int_{\R^n} |v|^2 (x,t)\, dx
\;\leq\; \frac{1}{2}&\int_{\R^n} |v|^2 (x,s)\,dx- \nu \int_s^t \int_{\R^n}
|\nabla v|^2 (x,\tau)\, dx\, d\tau\\
&\textrm{for almost every $s>0$ and for every $t>s$.}
\end{split}
\end{equation}
\end{theorem}

In what follows, the solutions of Theorem \ref{t:leray}
will be called {\em Leray solutions}. 
As is well known, Leray solutions are weakly continuous in time, i.e. 
\begin{equation}\label{e:weakcont}
t\;\mapsto\; \int_{\R^n} v (x,t)\cdot \varphi (x)\, dx 
\end{equation}
is continuous for every $\varphi\in L^2 (\R^n; \R^n)$. In other words
$v\in C ([0, T]; L^2_w (\R^n))$. More generally we have 

\begin{lemma}\label{l:weak}
Let $v$ be a weak solution of \eqref{euler} or a \emph{distributional} solution
of \eqref{navierstokes}, belonging to the space 
$L^\infty ([0, T]; L^2 (\R^n))$.
Then, $v$ can be redefined on a set of $t$ of measure zero
so that $v\in C ([0, T]; L^2_w (\R^n))$.
\end{lemma} 

This property (or a variant of it) is common to all distributional solutions of evolution equations
which can be written as balance laws (see for instance Theorem 4.1.1 in \cite{DafermosBook}) 
and can be proved by standard arguments. In Appendix A we include,
for the reader's convenience, a proof of a slightly more general statement,
which will be useful later. From now on we will use the
slightly shorter notation $C\bigl([0,T]; L^2_w\bigr)$ 
for $C\bigl([0,T]; L^2_w (\R^n)\bigr)$.

If a weak solution $v$ of \eqref{euler} 
is the strong limit of a sequence
of Leray solutions $v_k$ of \eqref{navierstokes} 
with vanishing viscosity $\nu = \nu_k\downarrow 0$,
then $v$ inherits in the limit \eqref{e:leray1} and \eqref{e:leray2}. 
Therefore one might say that this limit should be the weakest form of the energy inequality
that solutions of \eqref{euler} should satisfy.
This motivates the following definition.

\begin{definition}\label{d:weakenergy}
A weak solution $v\in C\bigl([0,T]; L^2_w \bigr)$ of \eqref{euler} satisfies the
{\em weak energy inequality} if
\begin{equation}\label{e:w1}
\int_{\R^n} |v|^2 (x,t)\, dx 
\;\leq\; \int_{\R^n} |v^0|^2 (x)\, dx\qquad\textrm{for every $t>0$,}
\end{equation}
and it satisfies the {\em strong energy inequality} if
\begin{equation}
\int_{\R^n} |v|^2 (x,t)\, dx 
\;\leq\; \int_{\R^n} |v|^2 (x,s)\, dx\qquad\textrm{for all $s,t$ with $t>s$.}
\end{equation}
Finally, $v$ satisfies the 
{\em energy equality}
if equality holds in \eqref{e:w1}.
\end{definition}

\subsection{The local energy inequality}\label{s:localenergy}
Consider next a Leray solution of \eqref{navierstokes}. 
Since $v\in L^\infty_t (L^2_x)$ and $\nabla v\in L^2_t (L^2_x)$,
the Sobolev inequality and a simple interpolation argument shows 
that $v\in L^3_{loc} (\R^n\times ]0,\infty[)$ 
if the space dimension $n$ is less or equal to $4$
{\footnote{Indeed, by the Sobolev embedding, we conclude that
$v\in L^2_t (L_x^{2^*})$. Interpolating between the spaces
$L^\infty L^2$ and $L^2L^{2^*}$ we conclude that $u\in L^r_t (L^s_x)$
for every exponents $r$ and $s$ satisfying the identities
$$
\frac{1}{r}=\frac{1-\alpha}{2} \quad \frac{1}{s} = \frac{\alpha}{2}
+ \frac{1-\alpha}{2^*} \;=\; \frac{1}{2} - \frac{1-\alpha}{n}
\qquad \mbox{for some $\alpha\in [0,1]$}.
$$
Plugging $\alpha = 2/(2+n)$ we obtain $r=s=2 (1+\frac{2}{n})=:q$.
Clearly, $q\geq 3$ for $n=2,3,4$.}}.
In this case, one could formulate a weak local form of the 
energy inequality, requiring that the natural
inequality corresponding to \eqref{LocEnNS}
holds in the distributional sense. This amounts to the condition 
\begin{equation}\label{locNS}
\int_0^{\infty}\int_{\R^n}|\nabla v|^2\varphi\,dx\,dt\;\leq\;\int_0^{\infty}\int_{\R^n} \frac{|v|^2}{2} (\partial_t \varphi + \nu \Delta \varphi)+ \left(\frac{|v|^2}{2} + p\right)v\cdot \nabla \varphi\,dx\,dt
\end{equation}
for any nonnegative $\varphi\in C^\infty_c (\R^n\times ]0,\infty[)$.
Note that, since $v\in L^3_{loc}$ 
and 
\begin{equation}\label{laplace}
\Delta p \;=\; \div \div (v\otimes v)\, ,
\end{equation}
by the Calderon--Zygmund estimates we have $p\in L^{3/2}_{loc}$.
Therefore $p v$ is a well--defined locally summable function.

It is not known whether the Leray solutions satisfy \eqref{locNS}. 
However, it is possible to construct global weak solutions satisfying
the weak energy inequality and the local energy inequality.
This fact has been proved for the first time by Scheffer in
\cite{Scheffer77} (see also the appendix of \cite{CKN}). 
The local energy
inequality is a fundamental ingredient in
the partial regularity theory initiated by Scheffer and
culminating in the work
of Caffarelli, Kohn and Nirenberg, see \cite{CKN} and \cite{Lin}.

\begin{theorem}\label{t:scheffer10}
Let $n\leq 4$ and let $v^0\in L^2 (\rn{n})$ 
be a divergence--free vector field.
Then there exists a weak solution $v$ of \eqref{navierstokes} 
with $\nabla v\in L^2_{loc}$ and which satisfies  
\eqref{e:leray1}, \eqref{e:leray2} and \eqref{locNS}.
\end{theorem}

By analogy, for weak solutions of \eqref{euler}, Duchon and 
Robert in \cite{DuchonRobert} have proposed to look at a local
form of the energy inequality \eqref{e:w1}.

\begin{definition}[Duchon--Robert]\label{d:DR} 
Consider an $L^3_{loc}$ 
weak solution $v$
of \eqref{euler}. We say that $v$ 
satisfies the {\em local energy inequality}
if 
\begin{equation}\label{e:loc1}
\partial_t \frac{|v|^2}{2} + \div \left(v \left( \frac{|v|^2}{2} + p\right)
\right)\;\leq\; 0
\end{equation}
in the sense of distributions, i.e. if
\begin{equation}\label{e:loc2}
\int_0^{\infty}\int_{\R^n} \frac{|v|^2}{2} \partial_t \varphi +\left(\frac{|v|^2}{2} + p\right)v\cdot \nabla \varphi
\;\geq\; 0
\end{equation}
for every nonnegative $\varphi\in C^\infty_c (\R^n\times ]0,\infty[)$.

Similarly, if the equality in \eqref{e:loc2} holds for every test function,
then we say that $v$ satisfies the {\em local energy equality}. 
\end{definition}

Since \eqref{laplace} holds even for weak solutions
of \eqref{euler}, $v\in L^3_{loc}$ implies
$p\in L^{3/2}_{loc}$, and hence the product $pv$
is well--defined.
Note, however, that, for solutions of
Euler, the requirement
$v\in L^3_{loc}$ is not at all natural, even in low dimensions: 
there is no apriori estimate yielding this property. 

\subsection{Measure--valued and dissipative solutions}\label{s:dissipative}

Two other very weak notions of solutions to incompressible Euler
have been proposed in the literature: DiPerna--Majda's
measure--valued solutions (see \cite{DipernaMajda1}) and
Lions' dissipative solutions (see Chapter 4.4 of \cite{Lions}). 

Both notions are based on considering weakly convergent sequences of Leray solutions of Navier-Stokes with vanishing viscosity.  

On the one hand, the possible oscillations in the nonlinear term $v\otimes v$ lead to the appearance of an additional term in the limit, where this term is subject to a certain pointwise convexity constraint. This can be formulated by saying that the weak limit is the barycenter of a measure--valued solution (cf.~\cite{DipernaMajda1} and also \cite{Bardos,TartarH} for alternative settings using Wigner- and H-measures). A closely related object is our "subsolution", defined in Section \ref{ss:functional}.

On the other hand, apart from the energy inequality, a version of the Gronwall inequality prevails in the weak limit, leading to the definition of dissipative solutions, cf.~Appendix B. As a consequence, 
dissipative solutions coincide with classical solutions as long as the latter exist:

\begin{theorem}[Proposition 4.1 in \cite{Lions}]\label{t:lions} 
If there exists a solution $v\in C ([0,T]; L^2 (\rn{n}))$
of \eqref{euler} such that $(\nabla v + \nabla v^T)\in L^1 ([0,T]; L^\infty
(\rn{n}))$, then any dissipative solution of \eqref{euler} is equal to
$v$ on $\rn{n}\times [0,T]$.
\end{theorem}

This is relevant for our discussion because of the following well known fact. 

\begin{proposition}\label{p:WImpDissBis}
Let $v\in C ([0, T]; L^2_w)$ be a weak solution of \eqref{euler}
satisfying the weak energy inequality. 
Then $v$ is a dissipative solution.
\end{proposition} 

Our construction yields initial data for which
the nonuniqueness results of Theorem \ref{t:main}
hold on any time interval $[0, \eps[$. However,
for sufficiently regular initial data, classical results
give the local existence of smooth solutions.  
Therefore, Proposition \ref{p:WImpDissBis} implies that,
{\em a fortiori}, the initial data considered
in our examples have necessarily a certain degree
of irregularity.

Though Proposition \ref{p:WImpDissBis} is well known,
we have not been able to find a 
reference for its proof and therefore
we include one in Appendix B (see the proof of
Proposition \ref{p:WImpDiss}).

\subsection{Admissible solutions to the $p$-system}
As usual, by a weak solution of \eqref{e:psistema} we understand
a pair $(\rho, v)\in L^{\infty}(\R^n)$ such that the following identities
hold for every test function $\psi,\varphi\in C^\infty_c (\rn{n}\times [0,\infty[)$:
\begin{equation}\label{e:test1}
\int_0^{\infty}\int_{\R^n} 
\Bigl[\rho \partial_t \psi + \rho v \cdot \nabla_x \psi\Bigr]\,dx\,dt\;+\;
\int_{\rn{n}} \rho^0 (x)\, \psi (x,0)\, dx=0,
\end{equation}
\begin{equation}\label{e:test2}
\int_0^{\infty}\int_{\R^n}
\Bigl[\rho v\cdot \partial_t\varphi + \rho \langle v\otimes v, \nabla \varphi\rangle\Bigr]\,dx\,dt
\;+\; \int_{\rn{n}} \rho^0 (x) v^0 (x)\cdot \varphi (x,0)\, dx\, =0.
\end{equation}
Admissible solutions have to satisfy an additional constraint.
Consider the internal energy $\eps: \rn{+}\to \rn{}$ given through
the law $p (r)= r^2 \eps' (r)$. Then admissible solutions of
\eqref{e:test1} have to satisfy the inequality
\begin{equation}\label{e:entropy}
\partial_t \left[\rho \eps (\rho) + \frac{\rho |v|^2}{2}\right]
+ {\rm div}_x \left[ \left(\rho \eps (\rho) + \frac{\rho |v|^2}{2} +
p (\rho)\right) v\right]\;\leq\; 0\, 
\end{equation}
in the sense of distributions (cf.~(3.3.18) and (3.3.21) of \cite{DafermosBook}).
More precisely

\begin{definition}\label{d:admissible} A weak solution of
\eqref{e:psistema} is admissible if the following inequality
holds for every nonnegative $\psi\in C^\infty_c (\rn{n}\times\rn{})$:
\begin{eqnarray}
&&\int_0^{\infty}\int_{\R^n} \left[\left(\rho \eps (\rho) + 
\frac{\rho |v|^2}{2}\right)
\partial_t\psi + \left(\rho \eps (\rho) + \frac{\rho |v|^2}{2} +
p (\rho)\right)v\cdot \nabla_x\psi\right]\nonumber\\
&+& \int_{\rn{n}} \left(\rho^0 \eps (\rho^0) + 
\frac{\rho^0 |v^0|^2}{2}\right)
\psi (\cdot,0)\;\geq\; 0\, .\label{e:entropy2}
\end{eqnarray}
\end{definition}

\section{A criterion for the existence of wild solutions}

In this section we state some criteria to recognize initial
data $v^0$ which allow for many weak solutions
of \eqref{euler} satisfying the weak, strong and/or local energy inequality.
In order to state it, we need to introduce some
of the notation already used in \cite{us}. 

\subsection{The Euler equation as a differential inclusion}
In particular, we state the following lemma
(compare with Lemma 2.1 of \cite{us}).
Here and in what follows we denote by
$\S^n$ the space of symmetric $n\times n$ matrices, by $\S^n_0$ the subspace of $\S^n$ of matrices with trace $0$, and by $I_n$ the $n\times n$ identity matrix.

\begin{lemma}\label{LRNC} 
Suppose $v\in L^2 (\R^n\times [0,T]; \R^n)$, 
$u\in  L^2 (\R^n\times [0,T]; \S^n_0)$, and 
$q$ is a distribution such that
\begin{equation}\label{LR}
\begin{split}
\partial_tv+\textrm{div }u+\nabla q&=0\, ,\\
\textrm{div }v&=0\, .
\end{split}
\end{equation}
If $(v,u,q)$ solve \eqref{LR} and in addition
\begin{equation}\label{NC}
u=v\otimes v-\frac{1}{n}|v|^2I_n\quad
\textrm{ a.e.~in }\R^n\times[0,T]\, ,
\end{equation}
then $v$ and $p:=q-\frac{1}{n}|v|^2$ 
solve \eqref{euler} distributionally. Conversely,
if $v$ and $p$ solve \eqref{euler} distributionally,
$v$, $u= v\otimes v - \frac{1}{n}|v|^2I_n$ and $q= p + \frac{1}{n}|v|^2$
solve \eqref{LR} and \eqref{NC}.
\end{lemma}

Next, for every $r\geq 0$, we consider the set of \emph{Euler states of speed $r$}
\begin{equation}\label{e:K}
K_r\;:=\;\left\{(v,u)\in\rn{n}\times\S^n_0:\,
u=v\otimes v-\frac{r^2}{n}I_n,\, |v|= r\right\}
\end{equation}
(cf.~Section of
\cite{us}, in particular (25) therein). 
Lemma \ref{LRNC} says simply that solutions to
the Euler equations can be viewed as evolutions on the manifold of Euler states
subject to the linear conservation laws \eqref{LR}.

Next, we denote by $K_r^{co}$ the convex hull
in $\rn{n}\times \S^n_0$ of $K_r$. 
In the following Lemma we give an explicit formula for $K_r^{co}$.
Since it will be often used in the sequel, we introduce the following notation.
For $v,w\in \rn{n}$ let $v\odot w$ denote the symmetrized tensor product, that is
\begin{equation}
v\odot w \;=\; \frac{1}{2} \bigl(v\otimes w + w\otimes v), 
\end{equation}
and let $v\trs w$ denote its traceless part, that is
\begin{equation}\label{e:tfstens}
v\trs w \;=\; \frac{1}{2} \bigl(v\otimes w + w\otimes v)-\frac{v\cdot w}{n}I_n. 
\end{equation}
Note that
$$
v\trs v \;=\; v\otimes v -\frac{|v|^2}{n} I_n
$$
and hence $K_r$ is simply
$$
K_r \;=\; \left\{ (v, v\trs v): |v|=r\right\}\, .
$$

\begin{lemma}\label{l:hull}
For any $w\in\S^n$ let $\lambda_{max}(w)$ denote
the largest eigenvalue of $w$. For $(v,u)\in \R^n\times\S^n_0$ let
\begin{equation}
e(v,u):=\frac{n}{2}\lambda_{max}(v\otimes v-u).
\end{equation}
Then
\begin{itemize}
\item[(i)] $e: \R^n\times\S^n_0\to\R$ is convex;
\item[(ii)] $\frac{1}{2}|v|^2\leq e(v,u)$, with equality if and only if $u=v\otimes v-\frac{|v|^2}{n}I_n$;
\item[(iii)] $|u|_{\infty}\leq 2\frac{n-1}{n}\, e(v,u)$, 
where $|u|_{\infty}$ denotes the operator norm of the matrix;
\item[(iv)] The $\frac{1}{2}r^2$--sublevel set of $e$ is the convex hull of
$K_r$, i.e.
\begin{equation}\label{e:hull}
K_r^{co} \;=\; \left\{ (v,u)\in \R^n\times \S^n_0 : 
e(v,u)\leq\frac{r^2}{2}\right\}\, .
\end{equation}
\item[(v)] If $(u,v)\in \R^n\times \S^n_0$, then $\sqrt{2e(v,u)}$
gives the smallest $\rho$ for which $(u,v)\in K_\rho^{co}$.
\end{itemize}
\end{lemma}

In view of (ii) if 
a triple $(v,u,q)$ solving \eqref{LR} corresponds a solution of the
Euler equations via the correspondence in Lemma \ref{LRNC},
then $e(v,u)$ is simply the energy density of the solution. In view of this
remark, if $(v,u,q)$ is a solution of \eqref{LR}, 
$e(v,u)$ will be called the {\em generalized energy density}, and 
$E(t) = \int_{\rn{n}} e (v(x,t), u(x,t)) dx$ will be called the {\em generalized energy}.

We postpone the proof of Lemma \ref{l:hull} to the next 
subsection and we state 
now the criterion for the existence of wild solutions. Its
proof, which is the core of the paper,
will be given in Section 4.

\begin{proposition}\label{p:main} 
Let $\Omega\subset \rn{n}$ be an open set (not necessarily bounded)
and let
$$
\bar{e}\in C 
\bigl(\overline{\Omega}\times ]0,T[\bigr)\cap 
C\bigl([0,T]; L^1(\Omega)\bigr).
$$
Assume there exists $(v_0,u_0,q_0)$ smooth solution of \eqref{LR} on 
$\R^n\times ]0,T[$ with the following properties:
\begin{equation}\label{e:(a)}
v_0\in C\bigl([0,T]; L^2_{w}\bigr),
\end{equation}
\begin{equation}\label{e:(b)}
\supp (v_0(\cdot,t), u_0 (\cdot,t))
\subset\subset\Omega\textrm{ for all }t\in]0,T[,
\end{equation}
\begin{equation}\label{e:(c)} 
e \bigl(v_0(x,t),u_0(x,t)\bigr)<\bar{e}(x,t)\textrm{ for all }(x,t)\in\Omega\times\, ]0,T[\,.
\end{equation}
Then there exist infinitely many weak solutions $v$ of the 
Euler equations \eqref{euler} in $\R^n\times [0,T[$ with pressure
\begin{equation}\label{e:(v)}
p=q_0-\frac{1}{n}|v|^2
\end{equation}
such that
\begin{eqnarray}
v&\in& C\bigl([0,T]; L^2_{w}\bigr),\label{e:(i)}\\
v(x,t)&=&v_0 (x, t)\label{e:(iii)}\quad\textrm{ for }t=0,T,\,\textrm{a.e. }x\in\R^n,\\
\frac{1}{2}|v(x, t)|^2&=&\bar{e}(x, t)\, {\mathbf 1}_{\Omega}\label{e:(iv)}\quad\textrm{ for every $t\in ]0,T[$, a.e.~$x\in\R^n$.}
\end{eqnarray}
\end{proposition}

\begin{remark}
The condition \eqref{e:(c)} implies that $\bar{e}>0$ on $\Omega\times]0,T[$. 
Hence $\overline{\Omega}\subset\R^n$ plays the role of the spatial support of the solutions. On the other hand, according to the statement of the Proposition 
the pair $(v,p)$ satisfies the Euler equations
\begin{eqnarray*}
\partial_tv+\div v\otimes v+\nabla p&=&0,\\
\div v&=&0,
\end{eqnarray*} 
 in all of $\R^n$ in the sense of distributions. In particular, even though the divergence-free condition
 implies that there is no jump of the normal trace of $v$ across the boundary $\partial\Omega$, the first equation shows that there is a jump of the normal trace of $v\otimes v$ which is compensated by a jump of $p$ across $\partial\Omega$.

\end{remark}
 
\subsection{Proof of Lemma \ref{l:hull}}

\begin{proof} {\bf (i)} Note that
\begin{eqnarray}
e(v,u) &=& \frac{n}{2}\max_{\xi\in S^{n-1}} \Bigl\langle\xi,(v\otimes v-u)\xi\Bigr\rangle
=\frac{n}{2}\max_{\xi\in S^{n-1}}\Bigl\langle\xi,\langle \xi,v\rangle v-u\xi\Bigr\rangle\nonumber\\
&=&\frac{n}{2}\max_{\xi\in S^{n-1}}\Bigl[|\langle \xi,v\rangle|^2-\langle \xi,u\xi\rangle\Bigr].
\label{e:equalitye}
\end{eqnarray}
Since for every $\xi\in S^{n-1}$ the map $(v,u)\mapsto |\langle \xi,v\rangle|^2-\langle\xi, u\xi\rangle$
is convex, it follows that $e$ is convex.

\medskip

{\bf (ii)} Since 
$v\otimes v=v\trs v+\frac{|v|^2}{n}I_n$,
we have, similarly to above, that
\begin{equation}\label{e:inequalitye}
\begin{split}
e(v,u)&=\frac{n}{2}\max_{\xi\in S^{n-1}} \Bigl\langle\xi,(v\trs v-u)\xi\Bigr\rangle+\frac{|v|^2}{2}\\
&=\frac{n}{2}\lambda_{max}(v\trs v-u)+\frac{|v|^2}{2}.
\end{split}
\end{equation}
Observe that, since $v\trs v-u$ is traceless, 
the sum of its eigenvalues is zero. Therefore 
$\lambda_{max}(v\trs v-u)\geq 0$ with equality if 
and only if $v\trs v-u=0$. 
This proves the claim.

\medskip

{\bf (iii)} From \eqref{e:equalitye} and \eqref{e:inequalitye} we deduce
\begin{equation*}
e(v,u)\geq \frac{n}{2}\max_{\xi\in S^{n-1}}\Bigl(-\langle \xi,u\xi\rangle\Bigr)= - \frac{n}{2} \lambda_{min} (u)\, . 
\end{equation*}
Therefore $-\lambda_{min}(u)\leq \frac{2}{n}e(v,u)$. 
Since $u$ is traceless, the sum of its eigenvalues is zero, hence 
\begin{equation*}
|u|_{\infty}\leq (n-1)|\lambda_{min}(u)|\leq \frac{2(n-1)}{n}\,e(v,u).
\end{equation*}

\medskip

{\bf (iv)} Without loss of generality we assume $r=1$. 
Let 
\begin{equation}\label{e:sublevel}
S_1\;:=\; \left\{(v,u)\in \R^n\times\S^n_0: e(v,u)\;\leq\; \frac{1}{2}\right\}\, .
\end{equation}
Observe that $e(v,u)=\frac{1}{2}$ for all $(v,u)\in K_1$, hence - by convexity of $e$ - 
$$
K_1^{co}\subset S_1.
$$
To prove the opposite inclusion, observe first of all that $S_1$ is convex by (i) and compact by (ii) and (iii). Therefore $S_1$ is equal to the closed convex hull of its extreme points. In light of this observation it suffices to show that the extreme points of $S_1$ are contained in $K_1$.

To this end let $(v,u)\in S_1\setminus K_1$. 
By a suitable rotation of the coordinate axes we may assume that $v\otimes v-u$ is diagonal, with diagonal entries $1/n\geq \lambda_1\geq\dots\geq\lambda_n$. Note that $(v,u)\notin K_1 \Longrightarrow \lambda_n<
1/n$. Indeed, if $\lambda_n =1/n$, then we have the identity
$u= v\otimes v - \frac{1}{n} I_n$. Since the trace of $u$ vanishes,
this identity implies $|v|^2=1$ and $u = v\otimes v - 
\frac{|v|^2}{n}I_n$, which give $(v,u)\in K_1$.

Let $e_1,\dots,e_n$ denote the coordinate unit vectors, and write 
$v=\sum_iv^ie_i$. Consider the pair $(\bar{v},\bar{u})\in 
\R^n\times\S^n_0$ defined by
\begin{equation*}
\bar{v}=e_n,\quad \bar{u}=\sum_{i=1}^{n-1}v^i(e_i\otimes e_n+e_n\otimes e_i).
\end{equation*}
A simple calculation shows that
$$
(v+t\bar{v})\otimes (v+t\bar{v})-(u+t\bar{u})=(v\otimes v-u)+(2t\,v^n+t^2)e_n\otimes e_n.
$$
In particular, since $\lambda_n<1/n$, $e(v+t\bar{v},u+t\bar{u})\leq 1/n$ for all sufficiently small $|t|$, so that $(v,u)+t(\bar{v},\bar{u})\in S_1$. This shows that $(v,u)$ cannot be an extreme point of $S_1$.

\medskip

{\bf (v)} is an easy direct consequence of (iv).
\end{proof} 

\section{Proof of Proposition \ref{p:main}}

Although the general strategy for proving
Proposition \ref{p:main} 
is based on Baire category arguments as in 
\cite{us}, there are several points in which Proposition 
\ref{p:main} differs, which give rise to 
technical difficulties. The main technical difficulty is 
given by the requirements \eqref{e:(i)} and \eqref{e:(iv)},
where we put a special emphasis on the fact that 
the equality in \eqref{e:(iv)} must hold
for {\it every} time $t$. The arguments in \cite{us}, which are 
based on the interplay between weak-strong convergence following 
\cite{Bernd}, yield only solutions in the space 
$L^{\infty}\bigl([0,T];L^2(\R^n)\bigr)$. Although such solutions 
can be redefined on a set of times of measure zero 
(see Lemma \ref{l:weak}) so that they belong to the space 
$C\bigl([0,T]; L^2_{w}\bigr)$, this gives the equality
\begin{equation}\label{e:almosteveryt}
\frac{1}{2}|v(\cdot, t)|^2=\bar{e}(\cdot, t)\, {\mathbf 1}_\Omega\textrm{ for almost every }t\in\,]0,T[\,.
\end{equation}
For the construction of 
solutions satisfying the strong energy inequality this
conclusion is not enough. Indeed, a consequence of 
Theorem \ref{t:main}c) is precisely
the fact that \eqref{e:(iv)} does not follow automatically
from \eqref{e:almosteveryt}.

This section is split into five parts. In 
\ref{ss:functional} we introduce the functional
framework, we state Lemma \ref{l:lsc}, Lemma 
\ref{l:funct=0} and Proposition \ref{p:perturbation},
and we show how Proposition \ref{p:main} follows from
them. The two lemmas are simple consequences of functional analytic
facts, and they are proved in \ref{ss:lemmas}.
Instead, the perturbation property of Proposition \ref{p:perturbation}
is the key point of the abstract argument, and it is the only place where the 
particular geometry of the equation enters. In 
\ref{ss:geometric} we introduce the waves which are the basic
building blocks for proving Proposition \ref{p:perturbation}. 
In \ref{ss:potential} we introduce a suitable
potential to localize the waves of \ref{ss:geometric}.
Finally, in \ref{ss:proofperturb}
we use these two tools and a careful construction
to prove Proposition \ref{p:perturbation}.

\subsection{Functional setup}\label{ss:functional}

We start by defining the space of "subsolutions" as follows. Let $v_0$ be a vectorfield as in Proposition \ref{p:main} with associated modified pressure $q_0$, and consider velocity fields 
$v:\R^n\times[0,T]\to\R^n$ which satisfy 
\begin{equation}\label{e:X0}
\div v=0, 
\end{equation}
the initial and boundary conditions
\begin{equation}\label{e:X1}
\begin{split}
v(x,0)&=v_0(x,0),\\
v(x,T)&=v_0(x,T),\\ 
\supp v(\cdot,t)&\subset\subset \Omega\textrm{ for all }t\in ]0,T[,
\end{split}
\end{equation}
and such that there exists a  
smooth matrix field $u:\R^n\times ]0,T[\to\S^n_0$ with 
\begin{equation}\label{e:X2}
\begin{split}
e\bigl(v(x,t),u(x,t)\bigr)<\bar{e}(x,t)\,&\textrm{ for all }(x,t)\in\Omega\times ]0,T[\,,\\
\supp u(\cdot,t)\subset\subset \Omega\, &\textrm{ for all } 
t\in ]0,T[\, ,\\
\partial_tv+\div u+\nabla q_0&=0\textrm{ in }\R^n\times [0,T].
\end{split}
\end{equation}

\begin{definition}[The space of subsolutions]\label{d:X_0}
Let $X_0$ be the set of such velocity fields, i.e.
$$
X_0=\Bigl\{v\in C^{\infty}\bigl(\R^n\times]0,T[\bigr)
\cap C\bigl([0,T];L^2_w\bigr):\,
\textrm{\eqref{e:X0},\eqref{e:X1},\eqref{e:X2} are satisfied}\Bigr\},
$$
and let $X$ be the closure of $X_0$ in $C\bigl([0,T]; L^2_{w}
\bigr)$. 
\end{definition}

We assume that $\bar{e}\in C\bigl([0,T];L^1(\Omega)\bigr)$, therefore there exists a constant $c_0$ such that
$\int_{\Omega}\bar{e}(x,t) dx\leq c_0$ for all $t\in[0,T]$. 
Since for any $v\in X_0$ we have
$$
\frac{1}{2}\int_{\R^n}|v(x,t)|^2\,dx\leq \int_{\Omega}\bar{e}(x,t) dx
\qquad \textrm{ for all }t\in [0,T],
$$
we see that $X_0$ consists of functions $v:[0,T]\to L^2(\R^n)$ 
taking values in a bounded 
subset $B$ of $L^2(\R^n)$. Without loss of generality
we can assume that $B$ is weakly closed.
Let $d_B$ be a metric on $B$ which metrizes the 
weak topology. Then $(B, d_B)$ 
is a compact metric space. Moreover, $d_B$ induces
naturally a metric
$d$ on the space $Y := C ([0,T]; (B, d_B))$ via
the definition
\begin{equation}\label{e:defd}
d (w_1, w_2)\;=\; \max_{t\in [0,T]} d_B (w_1 (\cdot, t),
w_2 (\cdot, t)).
\end{equation}
The topology induced by $d$ on $Y$ is equivalent 
to the topology of $Y$ as subset of
$C\bigl([0,T]; L^2_{w}\bigr)$. 
Moreover, by Arzel\`a-Ascoli's theorem, 
the space $(Y,d)$ is complete. Finally, 
$X$ is the closure in $(Y,d)$ of $X_0$, and hence $(X,d)$
is as well a complete metric space.

\begin{definition}[The functionals $I_{\eps, \Omega_0}$]
\label{d:functionals}
Next, for any $\eps>0$ and any bounded open set
$\Omega_0\subset \Omega$
consider the functional 
$$
I_{\eps, \Omega_0} (v):=\inf_{t\in [\eps,T-\eps]}\int_{\Omega_0}
\Bigl[\frac{1}{2}|v(x,t)|^2-\bar{e}(x,t)\Bigr]\,dx.
$$  
\end{definition}

It is clear that on $X$ each functional $I_{\eps,\Omega_0}$ 
is bounded from below. 

We are now ready to state the three important building
blocks of the proof of Proposition \ref{p:main}. 
The first two lemmas are simple consequences
of our functional analytic framework

\begin{lemma}\label{l:lsc}
The functionals $I_{\eps,\Omega_0}$ are lower-semicontinuous on $X$.
\end{lemma}

\begin{lemma}\label{l:funct=0}
For all $v\in X$ we have $I_{\eps, \Omega_0} (v)\leq 0$. If 
$I_{\eps, \Omega_0} (v)=0$ for every $\eps>0$ and every bounded
open set $\Omega_0\subset \Omega$, 
then $v$ is a weak solution of the Euler equations \eqref{euler} 
in $\R^n\times [0,T[$ with pressure 
$$
p=q_0-\frac{1}{n}|v|^2,
$$
and such that \eqref{e:(i)},\eqref{e:(iii)},\eqref{e:(iv)} are satisfied.
\end{lemma}

The following proposition is the key point in the whole argument, 
and it is the only place where the particularities of the equations 
enter. It corresponds to Lemma 4.6 of \cite{us}, though
its proof is considerably more complicated due to the special role
played by the time variable in this context.

\begin{proposition}[The perturbation property]\label{p:perturbation}
Let $\Omega_0$ and $\eps>0$ be given.
For all $\alpha>0$ there exists $\beta>0$ (possibly depending
on $\eps$ and $\Omega_0$) such that whenever $v\in X_0$ with
$$
I_{\eps, \Omega_0} (v)<-\alpha,
$$
there exists a sequence 
$v_k\in X_0$ with $v_k\overset{d}{\to} v$ and
$$
\liminf_{k\to\infty}I_{\eps, \Omega_0} (v_k)\geq I_{\eps, \Omega_0}
(v)+\beta\,. 
$$ 
\end{proposition}

\begin{remark}\label{rmk:perturbation}
In fact the proof of Proposition \ref{p:perturbation} will show that in case $\Omega$ is bounded and $\bar{e}$ is uniformly bounded in $\overline{\Omega}\times[0,T]$, the improvement $\beta$ in the statement can be chosen to be
$$
\beta=\min\{\alpha/2,C\alpha^2\},
$$
with $C$ only depending on $|\Omega|$ and $\Vert \bar{e}\Vert_{\infty}$.
\end{remark}

We postpone the proofs of these facts to
the following subsections, and now show how Proposition
\ref{p:main} follows from them and the general
Baire category argument. 

\begin{proof}[Proof of Proposition \ref{p:main}]
Since the functional $I_{\eps, \Omega_0}$ is lower-semicontinuous on 
the complete metric space $X$ and takes values in a bounded interval of
$\R$, it can be written as a pointwise supremum of
countably many continuous functionals, see 
Proposition 11 in Section 2.7 of Chapter IX of \cite{BourbakiGT}. 
Therefore, 
$I_{\eps, \Omega_0}$ is a Baire-1 map and hence its
points of continuity form a residual set in $X$. We claim
that if $v\in X$ is a point of continuity of $I_{\eps, \Omega_0}$, 
then $I_{\eps, \Omega_0} (v)=0$. 

To prove the claim, assume the contrary, i.e. that there exists $v\in X$ 
which is a point of continuity of $I_{\eps,\Omega_0}$ and 
$I_{\eps,\Omega_0}(v)<-\alpha$ for some $\alpha>0$. Choose a sequence 
$\{v_k\}\subset X_0$ such 
that $v_k\overset{d}{\to} v$. Then in particular $I_{\eps,\Omega_0}(v_k)\to I_{\eps,\Omega_0}(v)$ and so, by possibly renumbering the sequence, we may assume that $I_{\eps,\Omega_0}(v_k)<-\alpha$. Using Proposition \ref{p:perturbation} for each function $v_k$ and a standard diagonal argument, we find a new sequence $\{\tilde v_k\}\subset X_0$ 
such that
\begin{equation*}
\begin{split}
&\tilde v_k\overset{d}{\to} v\textrm{ in }X,\\
\lim_{k\to\infty}\,&I_{\eps,\Omega_0}(\tilde v_k)\geq I_{\eps,\Omega_0}(v)+\beta. 
\end{split}
\end{equation*}
This is in contradiction with the assumption that $v$ is a 
point of continuity of $I_{\eps,\Omega_0}$, thereby proving our claim.

Next, let $\Omega_k$ be an exhausting sequence of bounded open subsets
of $\Omega$. Consider the set $\Xi$ which is the intersection of 
$$
\Xi_k \;:=\; \left\{ v\in X\;:\; \mbox{$I_{1/k, \Omega_k}$ is continuous
at $v$}\right\}\, .
$$ 
$\Xi$ is the intersection of countably many residual sets and hence
it is residual. Moreover, if $v\in \Xi$, then $I_{\eps, \Omega_0} (v) =0$
for any $\eps>0$ and any bounded $\Omega_0\subset \Omega$. By Lemma 
\ref{l:funct=0}, any $v\in \Xi$ satisfies 
the requirements of Proposition \ref{p:main}. 
One can easily check that the cardinality of $X$ is infinite and therefore
the cardinality of any residual set in $X$ is infinite as well.
This concludes
the proof.
\end{proof}

\subsection{Proofs of Lemma \ref{l:lsc} and Lemma \ref{l:funct=0}} 
\label{ss:lemmas}

\begin{proof}[Proof of Lemma \ref{l:lsc}]
Assume for a contradiction that there exists $v_k,v\in X$ such that
$v_k\overset{d}{\to} v$ in $X$, but
\begin{eqnarray*}
&&\lim_{k\to\infty}\,
\inf_{t\in[\eps,T-\eps]}\int_{\Omega_0}
\Bigl[\frac{1}{2}|v_k(x,t)|^2-\bar{e}(x,t)\Bigr]dx\nonumber\\
&<& \inf_{t\in[\eps,T-\eps]}\int_{\Omega_0}\Bigl[\frac{1}{2}|v(x,t)|^2-\bar{e}(x,t)\Bigr]dx.
\end{eqnarray*}
Then there exists a sequence of times $t_k\in[\eps,T-\eps]$ such that
\begin{eqnarray}
&&\lim_{k\to\infty}\int_{\Omega_0}\Bigl[\frac{1}{2}|v_k(x,t_k)|^2-\bar{e}(x,t_k)\Bigr]\,dx\nonumber\\
&<& \inf_{t\in[\eps,T - \eps]}
\int_{\Omega_0}\Bigl[\frac{1}{2}|v(x,t)|^2-\bar{e}(x,t)\Bigr]\,dx.
\label{e:lsccontradict}
\end{eqnarray}
We may assume without loss of generality that $t_k\to t_0$. Since the convergence in $X$ is equivalent to the topology of $C\bigl([0,T]; 
L^2_{w}\bigr)$, we obtain that
$$
v_k(\cdot,t_k)\rightharpoonup v(\cdot,t_0)\,\textrm{ in }L^2(\R^n)\textrm{ weakly},
$$
and hence
$$
\liminf_{k\to\infty}\int_{\Omega_0}\Bigl[\frac{1}{2}|v_k(x,t_k)|^2-\bar{e}(x,t_k)\Bigr]\,dx\geq
\int_{\Omega_0}\Bigl[\frac{1}{2}|v(x,t_0)|^2-\bar{e}(x,t_0)\Bigr]\,dx.
$$
This contradicts \eqref{e:lsccontradict}, thereby concluding the proof.
\end{proof}

\begin{proof}[Proof of Lemma \ref{l:funct=0}]
For $v\in X_0$ there exists $u:\R^n\times ]0,T[\to\S^n_0$ 
such that \eqref{e:X2} holds. Therefore
$$
\frac{1}{2}|v(x,t)|^2\leq e\bigl(v(x,t),u(x,t)\bigr)<\bar{e}(x,t)
$$
for all $(x,t)\in\Omega\times ]0,T[$ and hence $I_{\eps, \Omega_0}
(v)\leq 0$ for $v\in X_0$. For general $v\in X$ 
the inequality follows from the density of $X_0$ and 
the lower-semicontinuity of $I_{\eps, \Omega_0}$.

Next, let $v\in X$ and assume that $I_{\eps, \Omega_0} (v)=0$
for every $\eps>0$ and every bounded open $\Omega_0\subset \Omega$. 
Let $\{v_k\}\subset X_0$ be a sequence such that 
$v_k\overset{d}{\to} v$ in $X$ and let $u_k$ be the associated 
sequence of matrix fields satisfying \eqref{e:X2}.  
The sequence $\{u_k\}$ satisfies the pointwise estimate 
\begin{equation*}
|u_k|_{\infty}\leq \frac{2(n-1)}{n}\,e(v_k,u_k)<
\frac{2(n-1)}{n}\,\bar{e}
\end{equation*}
in $\Omega$ because of Lemma \ref{l:hull} (iii), whereas $u_k=0$ outside $\Omega$. Therefore $\{u_k\}$ is locally uniformly bounded
in $L^\infty$ and hence, 
by extracting a weakly convergent  
subsequence and relabeling, we may assume that
$$
u_k\overset{*}{\rightharpoonup} u\,\textrm{ in }L^\infty_{loc}\bigl(\R^n\times]0,T[\bigr).
$$
Since $v_k\to v$ in $C\bigl([0,T];L^2_w\bigr)$ and 
$I_{\eps, \Omega_0} (v)=0$ for every choice of $\eps$ and $\Omega_0$, 
we see that $v$ satisfies \eqref{e:(i)}, 
\eqref{e:(iii)} and \eqref{e:(iv)}. Moreover, the linear equations
$$
\left\{\begin{aligned}
\partial_tv+\div u+\nabla q_0&=0,\\
\div v&=0
\end{aligned}\right.
$$
hold in the limit, and -- since $e$ is convex -- we have
\begin{equation}\label{e:limite}
e\bigl(v(x,t),u(x,t)\bigr) \leq \bar{e}(x,t)
\textrm{ for a.e. }(x,t)\in \Omega\times [0,T].
\end{equation}
To prove that $v$ is a weak solution of the Euler equations \eqref{euler} with pressure $p=q_0-\frac{1}{n}|v|^2$, in view of Lemma \ref{LRNC}, it suffices to show that 
\begin{equation}\label{e:exactsolution}
u=v\otimes v-\frac{|v|^2}{n}I_n\,\textrm{ a.e. in }\R^n\times [0,T].
\end{equation}
Combining \eqref{e:(iv)} and \eqref{e:limite} we have 
$$
\frac{1}{2}|v(x,t)|^2=e\bigl(v(x,t),u(x,t)\bigr)\textrm{ for almost every }(x,t)\in\Omega\times [0,T],
$$
so that \eqref{e:exactsolution} follows from
Lemma \ref{l:hull} (ii) and since $u=0$, $v=0$ outside $\Omega$.

\end{proof}

\subsection{Geometric setup}\label{ss:geometric}
In this subsection we introduce the first tool
for proving Proposition \ref{p:perturbation}.
The admissible segments defined below
correspond to suitable plane-wave solutions
of \eqref{LR}. More precisely, following L.~Tartar \cite{Tartar79},
the directions of these segments belong to 
the wave cone $\Lambda$ for 
the system of linear PDEs \eqref{LR}
(cf. Section 2 of \cite{us} and in particular
(7) therein).

\begin{definition}\label{def:admissible}
Given $r>0$ we will call $\sigma$ an \emph{admissible segment} if $\sigma$ is a line segment in $\R^n\times \S^n_0$ satisfying the following conditions:
\begin{itemize}
\item $\sigma$ is contained in the interior of $K_r^{co}$,
\item $\sigma$ is parallel to $(a,a\otimes a)-(b,b\otimes b)$ for some $a,b\in\R^n$ with $|a|=|b|=r$ and $b\neq \pm a$.
\end{itemize}
\end{definition} 

The following lemma, a simple consequence
of Carath\'eodory's theorem for convex sets, ensures the 
existence of sufficiently large admissible segments
(cf. with Lemma 4.3 of \cite{us}).

\begin{lemma}\label{l:geometric}
There exists a constant $C>0$, depending only on the dimension, such that for any $r>0$ and 
for any $(v,u)\in \textrm{int }K_r^{co}$ there exists an admissible line segment 
\begin{equation}\label{e:sigma}
\sigma=\Bigl[(v,u)-(\bar{v},\bar{u})\,,\,(v,u)+(\bar{v},\bar{u})\Bigr]
\end{equation}
such that
$$
|\bar{v}|\geq \frac{C}{r}(r^2-|v|^2).
$$
\end{lemma}

\begin{proof}
Let $z=(v,u)\in \textrm{int }K_r^{co}$. By Carath\'eodory's theorem $(v,u)$ lies in the interior of a simplex in $\R^n\times\S^n_0$ spanned by elements of $K_r$. In other words
$$
z=\sum_{i=1}^{N+1}\lambda_iz_i,
$$
where $\lambda_i\in\, ]0,1[\,$, $\sum_{i=1}^{N+1}\lambda_i=1$, 
$N=n(n+3)/2-1$ is the dimension of $\R^n\times\S^n_0$ and 
$$
z_i=\bigl(v_i,v_i\otimes v_i-\frac{r^2}{n}I_n\bigr)
$$
for some $v_i\in\R^n$ with $|v_i|=r$. By possibly perturbing the $z_i$ slightly, we can ensure that
$v_i\neq \pm v_j$ whenever $i\neq j$ (this is possible since $(v,u)$ is contained in the interior of the simplex).
Assume that the coefficients are ordered so that $\lambda_1=\max_i\lambda_i$. Then for any $j>1$
$$
z\,\pm\,\frac{1}{2}\lambda_j(z_j-z_1)\in \textrm{int }K_r^{co}.
$$
Indeed,
\begin{equation*}
z\,\pm\,\frac{1}{2}\lambda_j(z_j-z_1)=\sum_i\mu_iz_i,
\end{equation*}
where $\mu_1=\lambda_1\mp \frac{1}{2}\lambda_j$, $\mu_j=\lambda_j\pm\frac{1}{2}\lambda_j$ and $\mu_i=\lambda_i$ for $i\neq 1,j$. It is easy to see that $\mu_i\in\, ]0,1[$ for all $i=1\dots N$. 

On the other hand $z-z_1=\sum_{i=2}^{N+1}\lambda_i(z_i-z_1)$, so that
\begin{equation}\label{e:length}
|v-v_1|\leq N\max_{i=2\dots N+1}\lambda_i|v_i-v_1|.
\end{equation}
Let $j>1$ be such that $\lambda_j|v_j-v_1|=\max_{i=2\dots N+1}\lambda_i|v_i-v_1|$, and
let 
\begin{equation*}
\begin{split}
(\bar{v},\bar{u})&=\frac{1}{2}\lambda_j(z_j-z_1)\\
&=\frac{1}{2}\lambda_j\bigl(v_j-v_1,\,v_j\otimes v_j-v_1\otimes v_1\bigr).
\end{split}
\end{equation*}
Then $\sigma$, defined by \eqref{e:sigma}, is contained in the interior of $K_r^{co}$, hence it is an admissible segment. Moreover, by the choice of $j$ and using \eqref{e:length}
$$
\frac{1}{4rN}(r^2-|v|^2)= \frac{1}{4rN}(r+|v|)(r-|v|)\leq \frac{1}{2N}|v-v_1|\leq |\bar{v}|.
$$
This finishes the proof.
\end{proof}

\subsection{Oscillations at constant pressure}\label{ss:potential}

In this section we construct a potential for the linear conservation 
laws \eqref{LR}. Similar potentials were constructed in the
paper \cite{us} (see Lemma 3.4 therein). 
However, the additional feature
of this new potential is that it allows to
localize the oscillations at constant pressure,
which are needed in the proof of Proposition 
\ref{p:perturbation}.

As a preliminary step recall from Section 3 in \cite{us} that solutions of \eqref{LR} in $\R^n$ correspond to symmetric divergence--free matrix fields on $\R^{n+1}$ for which the $(n+1),(n+1)$ entry vanishes. To see this it suffices to consider the linear map
\begin{equation}\label{e:identifyU}
\R^n\times\S^n_0\times\R\ni(v, u, q) \quad\mapsto\quad
U=\begin{pmatrix} u+qI_n&v\\ v&0\end{pmatrix}.
\end{equation}
Note also that with this identification $q=\frac{1}{n}\tr U$. Therefore solutions of \eqref{LR} 
with $q\equiv 0$ correspond to matrix fields $U:\R^{n+1}\to\R^{(n+1)\times(n+1)}$ such that
\begin{equation}\label{e:potentialU}
\div U=0,\quad U^T=U,\quad U_{(n+1),(n+1)}=0,\quad \tr U=0.
\end{equation}
Furthermore, given a velocity vector $a\in\R^n$, the matrix of the corresponding Euler state is
$$
U_a=\begin{pmatrix}a\otimes a-\frac{|a|^2}{n}I_n&a\\ a&0\end{pmatrix}.
$$
The following proposition gives a potential for
solutions of \eqref{LR} oscillating
between two Euler states $U_a$ and $U_b$ of equal 
speed at constant pressure. 
 
\begin{proposition}\label{p:potential}
Let $a,b\in\R^n$ such that $|a|=|b|$ and $a\neq \pm b$. Then there exists a matrix--valued, constant coefficient, homogeneous linear differential operator of order 3
$$
A(\partial):C^{\infty}_c(\R^{n+1})\to C^{\infty}_c\bigl(\R^{n+1};\R^{(n+1)\times(n+1)}\bigr)
$$
and a space-time vector $\eta\in\R^{n+1}$ with the following properties:
\begin{itemize}
\item $U=A(\partial)\phi$ satisfies \eqref{e:potentialU} for all $\phi\in C^{\infty}_c(\R^{n+1})$
\item $\eta$ is not parallel to $e_{n+1}$;
\item if $\phi(y)=\psi(y\cdot\eta)$, then
$$
A(\partial)\phi(y)=(U_a-U_b)\,\psi'''(y\cdot\eta).
$$
\end{itemize}
\end{proposition}

\begin{proof}
A matrix valued homogeneous polynomial of degree 3
$$
A:\R^{n+1}\to\R^{(n+1)\times (n+1)}
$$
gives rise to a differential operator required by the proposition if and only if $A=A(\xi)$ satisfies
\begin{equation}\label{e:potentialA}
A\xi=0,\quad A^T=A,\quad Ae_{(n+1)}\cdot e_{(n+1)}=0,\quad \tr A=0
\end{equation}
for all $\xi\in\R^{n+1}$.

Define the $(n+1)\times (n+1)$ antisymmetric matrices
\begin{equation*}
\begin{split}
R&=a\otimes b - b\otimes a,\\
Q(\xi)&=\xi\otimes e_{n+1}-e_{n+1}\otimes \xi,
\end{split}
\end{equation*}
where in the definition of $R$ we treat $a,b\in\R^n$ as elements of $\R^{n+1}$ by setting the $(n+1)$'s coordinate zero.
The following facts are easily verified:
\begin{enumerate}[(i)]
\item $R\xi\cdot \xi=0$,\quad $Q(\xi)\xi\cdot \xi=0$, due to antisymmetry;
\item $R\xi\cdot e_{n+1}=0$, since $a\cdot e_{n+1}=b\cdot e_{n+1}=0$;
\item $R\xi\cdot Q(\xi)\xi = 0$, because by (i) and (ii) $R\xi$ 
is perpendicular to the range of $Q$.
\end{enumerate}
Let
\begin{equation*}
A(\xi)=R\xi\odot \bigl(Q(\xi)\xi\bigr)=\frac{1}{2}\Bigl(R\xi\otimes \bigl(Q(\xi)\xi\bigr)+\bigl(Q(\xi)\xi\bigr)\otimes R\xi   \Bigr)
\end{equation*}
The properties (i),(ii),(iii) immediately imply \eqref{e:potentialA}. 

Now define $\eta\in\R^{n+1}$ by
$$
\eta=\frac{-1}{(|a||b|+a\cdot b)^{2/3}}\biggl(a+b-(|a||b|+a\cdot b)e_{n+1}\biggr).
$$
Since $|a|=|b|$ and $a\neq\pm b$, $|a||b|+a\cdot b\neq 0$ so that $\eta$ is well--defined and non--zero.
Moreover, a direct calculation shows that
$$
A(\eta)=\begin{pmatrix}a\otimes a-b\otimes b& a-b\\ a-b& 0\end{pmatrix}=U_a-U_b.
$$
Finally, observe that if $\phi(y)=\psi(y\cdot\eta)$, then
$A(\partial)\phi(y)=A(\eta)\psi'''(y\cdot\eta)$.
\end{proof}

The following simple lemma ensures that the oscillations of the plane-waves 
produced by Proposition \ref{p:potential} have a certain size in terms of functionals 
of the type $I_{\eps,\Omega_0}$. 

\begin{lemma}\label{l:trigonometry}
Let $\eta\in\R^{n+1}$ be a vector which is not parallel to $e_{n+1}$. Then for any
bounded open set $B\subset\R^n$
$$
\lim_{N\to\infty}\int_B\sin^2\bigl(N\eta\cdot(x,t)\bigr)\,dx\;=\;\frac{1}{2}|B|
$$
uniformly in $t\in\R$.
\end{lemma}

\begin{proof}
Let us write $\eta=(\eta',\eta_{n+1})\in\R^n\times\R$, so that $\eta'\in\R^n\setminus\{0\}$. By elementary trigonometric identities
\begin{equation*}
\begin{split}
\sin^2\bigl(N\eta&\cdot(x,t)\bigr)=
\sin^2(N\eta'\cdot x)+\\
&+\sin^2(N\eta_{n+1}t)\cos(2N\eta'\cdot x)+\frac{1}{2}\sin(2N\eta'\cdot x)\sin(2N\eta_{n+1}t).
\end{split}
\end{equation*}
For the second term we have
\begin{equation*}
\Bigl|\int_B\sin^2(N\eta_{n+1}t)\cos(2N\eta'\cdot x)dx\Bigr|
\leq \Bigl|\int_{B}\cos(2N\eta'\cdot x)dx\Bigr|\,\to 0
\end{equation*}
as $N\to \infty$, and similarly the third term vanishes in the limit uniformly in $t$. The statement of the lemma now follows easily. 
\end{proof}

%%%%%%%%%%%%%%%%%%%%%%%%%%%%%%%%%%%%%%

\subsection{Proof of the perturbation property}\label{ss:proofperturb}
We are now ready to conclude the proof of Proposition
\ref{p:perturbation}.

\medskip

\begin{figure}[htbp]
\begin{center}
    \input{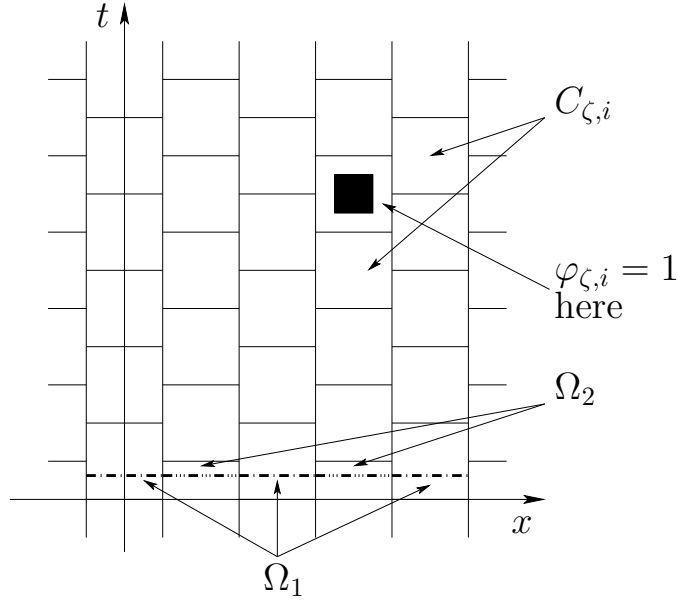}
    \caption{The ``shifted'' grid in dimension $1+1$.}
    \label{f:grid}
\end{center}
\end{figure}

{\bf Step 1. Shifted grid. } 
We start by defining a grid on $\R^n_x\times\R_t$ of size $h$. For $\zeta\in\Z^n$ 
let $|\zeta|=\zeta_1+\dots +\zeta_n$ and let $Q_{\zeta},\tilde Q_{\zeta}$ be cubes in $\R^n$ centered at $\zeta h$ with sidelength $h$ and $\frac{3}{4}h$ respectively, i.e.
$$
Q_{\zeta}:=\zeta h + \left[-\frac{h}{2},\frac{h}{2}\right]^n,\;
\tilde Q_{\zeta}:=\zeta h + \left[-\frac{3h}{8},\frac{3h}{8}\right]^n\,.
$$
Furthermore,
for every $(\zeta,i)\in\Z^n\times\Z$ let
$$
C_{\zeta,i}=\begin{cases} Q_{\zeta}\times [ih,(i+1)h]&
\textrm{ if $|\zeta|$ is even},\\
Q_{\zeta}\times [(i-\frac{1}{2})h,(i+\frac{1}{2})h]&
\textrm{ if $|\zeta|$ is odd.}\end{cases}
$$
Next, we let $0\leq\varphi\leq 1$ 
be a smooth cutoff function on $\R^n_x\times\R_t$, with support contained
in $[-h/2, h/2]^{n+1}$,
identically $1$ on $[-3h/8, 3h/8]^{n+1}$ and strictly less than $1$ outside. 
Denote by $\varphi_{\zeta, i}$ the obvious
translation of $\varphi$ supported in $C_{\zeta, i}$, and let
$$
\phi^h \;:=\; \sum_{\zeta\in \Z^n, i\in \Z} \varphi_{\zeta, i}\, .
$$
Given an open and bounded set $\Omega_0$, let
\begin{equation*}
\Omega^h_1=\bigcup\bigl\{\tilde Q_{\zeta}:\,|\zeta|\textrm{ even, }Q_{\zeta}\subset\Omega_0\bigr\}
,\quad \Omega^h_2=\bigcup\bigl\{\tilde Q_{\zeta}:\,|\zeta|\textrm{ odd, }Q_{\zeta}\subset\Omega_0
\bigr\}\, .
\end{equation*}
Observe that
$$
\lim_{h\to 0}|\Omega^h_{\nu}|\;=\; \frac{1}{2}\left(\frac{3}{4}\right)^n|\Omega_0|\quad\textrm{ for }\nu=1,2,
$$ 
and for every fixed $t$ the set
$\{x\in \Omega_0: \phi^h (x,t)=1\}$ contains at least one of the sets $\Omega^h_{\nu}$, see Figure \ref{f:grid}. Indeed, if
$$
\tau^h_1=\bigcup_{i\in\N}\Bigl[(i+\frac{1}{4})h,(i+\frac{3}{4})h\Bigr[\;\textrm{ and }\,
\tau^h_2=\bigcup_{i\in\N}\Bigl[(i-\frac{1}{4})h,(i+\frac{1}{4})h\Bigr[\,,
$$
then $\tau^h_1\cup\tau^h_2=\R$, and for $\nu=1,2$
$$
\phi^h(x,t)=1\textrm{ for all }(x,t)\in\Omega^h_{\nu}\times \tau^h_{\nu}.
$$

\medskip

Now let $v\in X_0$ with
$$
I_{\eps,\Omega_0}(v)<-\alpha
$$
for some $\alpha>0$, and let $u:\Omega\times\,]0,T[\,\to\S^n_0$ be a corresponding smooth matrix field satisfying \eqref{e:X2}. Let
\begin{equation}\label{e:defM}
M=\max_{\Omega_0\times[\eps/2,T-\eps/2]} \bar{e},
\end{equation}
and let $E_h:\Omega_0\times[\eps,T-\eps]\to\R$ be 
the step-function on the grid defined by
$$
E_h(x,t)=E_h(\zeta h, ih)=\frac{1}{2}\bigl|v(\zeta h,ih)\bigr|^2-\bar{e}(\zeta h,ih)\quad\textrm{ for }(x,t)\in C_{\zeta,i}.
$$
This is well--defined provided $h<\eps$.
Since $v$ and $\bar{e}$ are uniformly continuous on $\Omega_0\times [\eps/2,T-\eps/2]$, 
for any $\nu\in\{1,2\}$
\begin{equation*}
\lim_{h\to 0}\int_{\Omega^h_{\nu}}E_h(x,t) dx\;=\;\frac{1}{2}\left(\frac{3}{4}\right)^n
\int_{\Omega_0} \Bigl[\frac{1}{2}|v(x,t)|^2-\bar{e}(x,t)\Bigr]dx
\end{equation*}
uniformly in $t\in [\eps, T-\eps]$.
In particular there exists a dimensional constant $c>0$ such that, for all sufficiently small grid sizes $h$ and for any 
$t\in [\eps,T-\eps]$, we have
\begin{equation}\label{e:needed}
\begin{split}
\int_{\Omega^h_{\nu}}&|E_h(x,t)| dx\geq c\alpha\\
&\textrm{ whenever }\int_{\Omega_0}\Bigl[\frac{1}{2}|v(x,t)|^2-\bar{e}(x,t)\Bigr]dx\leq -\frac{\alpha}{2}.
\end{split}
\end{equation}
 
\medskip

Next, for each $(\zeta,i)\in\Z^n\times\Z$ such that $C_{\zeta, i}\subset\Omega_0\times 
[\eps/2, T-\eps/2]$ let 
$$
z_{\zeta,i}=\bigl( v(\zeta h,ih),u(\zeta h,ih)\bigr),
$$
and, using Lemma \ref{l:geometric}, choose a segment 
$$
\sigma_{\zeta, i}=\bigl[ z_{\zeta,i}-\bar{z}_{\zeta,i},z_{\zeta,i}+\bar{z}_{\zeta,i}\bigr]
$$
admissible for $r= \sqrt{2 \bar{e} (\zeta h,ih)}$ (cf.~Definition \ref{def:admissible}) with midpoint
$z_{\zeta,i}$ and direction 
$\bar{z}_{\zeta,i}=\bigl(\bar{v}_{\zeta,i},\bar{u}_{\zeta,i}\bigr)$ such that
\begin{equation}\label{e:length2}
|\bar{v}_{\zeta,i}|^2\geq \frac{C}{\bar{e}(\zeta h,ih)}|E_h(\zeta h,ih)|^2\geq 
\frac{C}{M}|E_h(\zeta h,ih)|^2.
\end{equation}
Since $z:=(v,u)$ and $\bar{e}$ are uniformly continuous, for sufficiently small $h$ we have
\begin{equation}\label{e:suffsmallh}
e\bigl( z(x,t)+\lambda \bar{z}_{\zeta,i}\bigr)<\bar{e}(x,t)\quad\textrm{ for all }\lambda\in [-1,1]\textrm{ and }(x,t)\in C_{\zeta,i}.
\end{equation}
Thus we fix the grid size $0<h<\eps/2$ so that the estimates \eqref{e:needed} and \eqref{e:suffsmallh} hold.

\bigskip

{\bf Step 2. The perturbation. } 
Fix $(\zeta, i)$ for the moment. Corresponding to the admissible segment 
$\sigma_{\zeta,i}$,
in view of Proposition \ref{p:potential} and the identification \eqref{e:identifyU} there exists an operator $A_{\zeta,i}$ and a direction $\eta_{\zeta,i}\in\R^{n+1}$, not parallel to $e_{n+1}$, such that for any $N\in\N$
$$
A_{\zeta,i} \left(N^{-3}\cos \left(N\eta_{\zeta,i}\cdot (x,t)\right)\right)
\;=\; \bar{z}_{\zeta,i} \sin \left(N \eta_{\zeta,i}\cdot (x,t)\right)\,,
$$
and such that the pair $(v_{\zeta,i}, u_{\zeta, i})$ defined by
$$
(v_{\zeta,i}, u_{\zeta, i}) (x,t)\;:=\;A_{\zeta,i} \Bigl[\varphi_{\zeta,i}(x,t)\, N^{-3}\cos \left(N \eta_{\zeta,i}\cdot (x,t)\right)\Bigr]
$$
satisfies \eqref{LR} with $q\equiv 0$. Note that
$(v_{\zeta, i},u_{\zeta,i})$ is supported in the cylinder 
$C_{\zeta, i}$ and that
\begin{equation}\label{e:linfty}
\begin{split}
&\Bigl\|(v_{\zeta, i},u_{\zeta,i}) - \varphi_{\zeta, i} \bar{z}_{\zeta,i} \sin \left(N \eta_{\zeta,i}\cdot (x,t)\right)\Bigr\|_\infty\\
=\; &\Bigl\|A_{\zeta,i} \Bigl[\varphi_{\zeta,i}\, N^{-3}\cos \left(N \eta_{\zeta,i}\cdot (x,t)\right)\Bigr]\\
&\qquad-\; \varphi_{\zeta,i}\,A_{\zeta,i} \Bigl[N^{-3}\cos \left(N \eta_{\zeta,i}\cdot (x,t)\right)\Bigr]\Bigr\|_{\infty}\\
\leq\; & C\bigl(A_{\zeta,i},\eta_{\zeta,i},\|\varphi_{\zeta,i}\|_{C^3}\bigr)\frac{1}{N},
\end{split}
\end{equation}
since $A_{\zeta,i}$ is a linear differential operator 
of homogeneous degree 3. Let
$$
(\tilde v_N,\tilde u_N)\;:=\;
\sum_{(\zeta,i):C_{\zeta,i}\subset\Omega_0\times [\eps,T-\eps]}(v_{\zeta,i},u_{\zeta,i})
$$
and
$$
(v_N,u_N)=(v,u)+ (\tilde{v}_N, \tilde{u}_N)\, .
$$
Observe that the sum consists of finitely many terms. Therefore from
\eqref{e:suffsmallh} and \eqref{e:linfty} we deduce that
there exists $N_0\in\N$ such that
\begin{equation}\label{e:inX_0}
v_N\in X_0\textrm{ for all }N\geq N_0.
\end{equation}
Furthermore, recall that for all $(x,t)\in \Omega_{\nu}\times\tau_{\nu}$ we have $\phi^h(x,t)=1$ and hence
$$
|\tilde v_N(x,t)|^2=|\bar{v}_{\zeta,i}|^2\sin^2(N\eta_{\zeta,i}\cdot (x,t))\,,
$$ 
where $i\in\N$ is determined by the inclusion $(x,t)\in C_{\zeta,i}$. 
Since $\eta_{\zeta,i}\in\R^{n+1}$ is not parallel to $e_{n+1}$, from Lemma \ref{l:trigonometry} 
we see that
$$
\lim_{N\to\infty}\int_{\tilde Q_{\zeta}}|\tilde v_N(x,t)|^2dx\,=\,\frac{1}{2}\int_{\tilde Q_{\zeta}}|\bar{v}_{\zeta,i}|^2dx
$$
uniformly in $t$. In particular, using \eqref{e:length2} and summing over all $(\zeta,i)$ such that $C_{\zeta,i}\subset \Omega_0\times[\eps,T-\eps]$, we obtain
\begin{equation}\label{e:size}
\lim_{N\to\infty}\int_{\Omega^h_{\nu}}\frac{1}{2}|\tilde v_N(x,t)|^2dx\geq \frac{c}{M}\int_{\Omega_{\nu}^h}|E_h(x,t)|^2dx
\end{equation}
uniformly in $t\in\tau_{\nu}\cap[\eps,T-\eps]$,
where $c>0$ is a dimensional constant.

\bigskip

{\bf Step 3. Conclusion.} 
For each $t\in [\eps,T-\eps]$ we have
\begin{equation*}
\begin{split}
\int_{\Omega_0}\Bigl[\frac{1}{2}&|v_N(x,t)|^2-\bar{e}(x,t)\Bigr]dx=
\int_{\Omega_0}\Bigl[\frac{1}{2}|v(x,t)|^2-\bar{e}(x,t)\Bigr]dx\\
&+\int_{\Omega_0}\frac{1}{2}|\tilde v_N(x,t)|^2dx+\int_{\Omega_0}\tilde v_N(x,t)\cdot v(x,t)dx.
\end{split}
\end{equation*}
Since $v$ is smooth on $\Omega_0\times [\eps/2,T-\eps/2]$, 
\begin{equation*}
\int_{\Omega_0}\tilde v_N(x,t)\cdot v(x,t)dx\;\to\;0\textrm{ as }N\to\infty,\textrm{ uniformly in $t$},
\end{equation*}
hence
\begin{equation*}
\liminf_{N\to\infty}I_{\eps,\Omega_0}(v_N)\geq\liminf_{N\to\infty}\inf_{t\in[\eps,T-\eps]}
\Biggl\{\int_{\Omega_0}\Bigl[\frac{1}{2}|v|^2-\bar{e}\Bigr]dx+\int_{\Omega_0}\frac{1}{2}|\tilde v_N|^2dx\Biggr\}.
\end{equation*}
Since the limit in \eqref{e:size} is uniform in $t$, it follows that
\begin{equation*}
\begin{split}
&\liminf_{N\to\infty}I_{\eps,\Omega_0}(v_N)\geq  \inf_{t\in[\eps,T-\eps]}\Biggl\{\int_{\Omega_0}\Bigl[\frac{1}{2}|v|^2-\bar{e}\Bigr]dx+\frac{c}{M}\min_{\nu\in\{1,2\}}\int_{\Omega_{\nu}^h}|E_h|^2dx\Biggr\}\\
&\quad\geq  \inf_{t\in[\eps,T-\eps]}\Biggl\{\int_{\Omega_0}\Bigl[\frac{1}{2}|v|^2-\bar{e}
\Bigr]dx+\frac{c}{M|\Omega_0|}\min_{\nu\in\{1,2\}}\Bigl(\int_{\Omega_{\nu}^h}|E_h|dx\Bigr)^2\Biggr\},
\end{split}
\end{equation*}
where we have applied the Cauchy-Schwarz inequality on the last integral.
We conclude, using \eqref{e:needed}, that
\begin{equation*}
\begin{split}
\liminf_{N\to\infty}I_{\eps,\Omega_0}(v_N)&\geq \min\biggl\{-\frac{\alpha}{2},-\alpha+\frac{c}{M|\Omega_0|}\alpha^2\biggr\}\\
&\geq -\alpha+\min\biggl\{\frac{\alpha}{2},\frac{c}{M|\Omega_0|}\alpha^2\biggr\}.
\end{split}
\end{equation*}
On the other hand we recall from \eqref{e:inX_0} that $v_N\in X_0$ for $N\geq N_0$ and furthermore clearly $v_N\overset{d}{\to} v$. This concludes the proof.

%%%%%%%%%%%%%%%%%%%%%%%%%%%%%%%%%%%%%%%%%%

\section{Construction of suitable initial data}

In this section we construct examples of initial data
for which we have a ``subsolution'' in the sense of Proposition
\ref{p:main}. We fix here a
bounded open set $\Omega\subset\R^n$.

\begin{proposition}\label{p:initialdata} There exist  
triples $(\bar{v},\bar{u},\bar{q})$ solving 
\eqref{LR} in $\R^n\times\rn{}$ and enjoying the 
following properties:
\begin{equation}\label{e:init(a)}
\bar{q}\equiv 0,\,  (\bar{v}, \bar{u})
\mbox{ is smooth in $\rn{n}\times (\rn{}\setminus \{0\})$ and }
\bar{v}\in C\bigl(\rn{}; L^2_{w}\bigr)\, ,
\end{equation}
\begin{equation}\label{e:init(b')}
\supp (\bar{v} ,\bar{u})\subset\overline{\Omega}\times\,]-T,T[\,,
\end{equation}
\begin{equation}\label{e:init(b)}
\supp (\bar{v} (\cdot,t),\bar{u} (\cdot, t))
\subset\subset\Omega\textrm{ for all }t\neq 0\,,
\end{equation}
\begin{equation}\label{e:init(c)} 
e \bigl(\bar{v}(x,t),\bar{u} (x,t)\bigr)<1
\textrm{ for all }(x,t)\in \rn{n}\times (\rn{}\setminus\{0\})\,.
\end{equation}
Moreover
\begin{equation}\label{e:init(d)}
\frac{1}{2}|\bar{v} (x,0)|^2=1\textrm{ a.e. in }\Omega.
\end{equation}
\end{proposition}

\begin{remark}
Observe that \eqref{e:init(c)} and \eqref{e:init(d)} together imply that $\bar{v}(t)\to \bar{v}(0)$ {\em strongly} in $L^2(\R^n)$ as $t\to 0$.
\end{remark}

\begin{proof}
In analogy with Definition \ref{d:X_0} we consider the space $X_0$, defined as the set of vector fields $v:\R^n\times ]-T,T[\to \R^n$ in $C^{\infty}(\R^n\times ]-T,T[)$ to which there exists 
a smooth matrix field $u:\R^n\times ]-T,T[\to\S^n_0$ such that
\begin{equation}\label{e:id-1}
\begin{split}
\div v&=0,\\
\partial_t v+\div u&=0,
\end{split}
\end{equation}
\begin{equation}\label{e:id-2}
\supp (v,u)\subset\Omega\times [-T/2,T/2[\,,
\end{equation}
and
\begin{equation}\label{e:id-3}
e\bigl(v(x,t),u(x,t)\bigr)<1\quad\textrm{ for all }(x,t)\in\Omega\times ]-T,T[\,.
\end{equation}
This choice of $X_0$ corresponds - up to changing the time interval under consideration - in Section \ref{ss:functional} to the choices $(v_0,u_0,q_0)\equiv (0,0,0)$ 
and $\bar{e}\equiv 1$. Similarly to before, $X_0$ consists of functions 
$v:]-T,T[\,\to L^2(\R^n)$ taking values in a bounded set $B\subset L^2(\R^n)$ (recall that in this section we assume $\Omega$ is bounded). On $B$ the weak topology of $L^2$ is metrizable, and correspondingly we find a metric $d$ on $C(]-T,T[\,,B)$ inducing the topology
of $C(]-T,T[\,,L^2_w(\R^n))$. 

Next we note that with minor modifications the proof of the perturbation property in Section \ref{ss:proofperturb} leads to the following claim (cf.~Remark \ref{rmk:perturbation} following the statement of Proposition \ref{p:perturbation}):

\bigskip

\noindent{\bf Claim: }
Let $\Omega_0\subset\subset\Omega$ be given. Let $v\in X_0$ with associated matrix field $u$ and let $\alpha>0$ such that
$$
\int_{\Omega_0}\Bigl[\frac{1}{2}|v(x,0)|^2-1\Bigr]\,dx<-\alpha.
$$
Then for any $\eps>0$ there exists a sequence $v_k\in X_0$ with associated smooth matrix field $u_k$ such that
\begin{equation}
\supp(v_k-v,u_k-u)\subset \Omega_0\times [-\eps,\eps],
\end{equation}
\begin{equation}\label{e:id-weak}
v_k\overset{d}{\to} v,
\end{equation}
and
\begin{equation}
\liminf_{k\to\infty}\int_{\Omega_0}\frac{1}{2}|v_k(x,0)|^2\,dx\geq \int_{\Omega_0}\frac{1}{2}|v(x,0)|^2\,dx+\min\biggl\{\frac{\alpha}{2},C\alpha^2\biggr\},
\end{equation}
where $C$ is a fixed constant independent of $\eps,\alpha,\Omega_0$ and $v$.

\bigskip

Fix an exhausting sequence of bounded open subsets $\Omega_k\subset\Omega_{k+1}\subset\Omega$, each compactly contained in $\Omega$, and such that $|\Omega_{k+1}\setminus\Omega_k|\leq 2^{-k}$. Let also $\rho_{\eps}$ be a standard mollifying kernel in $\R^n$. Using the claim above we construct inductively a sequence of velocity fields $v_k\in X_0$, associated matrix fields $u_k$ and a sequence of numbers 
$\eta_k<2^{-k}$ as follows. 

First of all let $v_1\equiv 0$ and $u_1\equiv 0$.  
Having obtained $(v_1,u_1),\dots,(v_k,u_k)$ and $\eta_1,\dots,\eta_{k-1}$ we choose 
$\eta_k<2^{-k}$ in such a way that
\begin{equation}\label{e:strongmolly}
\Vert v_k-v_k*\rho_{\eta_k}\Vert_{L^1}<2^{-k}.
\end{equation}
Furthermore, we define
$$
\alpha_k=-\int_{\Omega_k}\Bigl[\frac{1}{2}|v_k(x,0)|^2-1\Bigr]\,dx.
$$
Note that due to \eqref{e:id-3} we have $\alpha_k>0$.

Then we apply the claim with $\Omega_k$, $\alpha=\frac{3}{4}\alpha_k$  and 
$\eps=2^{-k}T$ to obtain 
$v_{k+1}\in X_0$ and associated smooth matrix field $u_{k+1}$ such that
\begin{equation}\label{e:id-supp}
\supp (v_{k+1}-v_k,u_{k+1}-u_k)\subset\Omega_k\times \bigl[-2^{-k}T,2^{-k}T\bigr]\,,
\end{equation}
\begin{equation}\label{e:id-conv}
d(v_{k+1},v_k)<2^{-k}\,,
\end{equation}
\begin{equation}\label{e:id-norm}
\int_{\Omega_k}\frac{1}{2}|v_{k+1}(x,0)|^2dx\geq \int_{\Omega_k}\frac{1}{2}|v_k(x,0)|^2dx+\frac{1}{4}\min\{\alpha_k,C\alpha_k^2\},
\end{equation}
and recalling that $d$ induces the topology of $C(]-T,T[,L^2_w)$ we can prescribe in addition that
\begin{equation}\label{e:id-weak2}
\Vert (v_k-v_{k+1})*\rho_{\eta_j}\Vert_{L^2(\Omega)}<2^{-k}\textrm{ for all }j\leq k\textrm{ for }t=0.
\end{equation}

>From \eqref{e:id-conv} we deduce that there exists $\bar{v}\in C(]-T,T[\,,L^2_w(\Omega))$ such that
$$
v_k\overset{d}{\to}\bar{v}.
$$
>From \eqref{e:id-supp} we see that for any compact subset of 
$\Omega\times\,]-T, 0[\,\cup\, ]0,T[$ there exists $k_0$ such that $(v_k,u_k)=(v_{k_0},u_{k_0})$ for all $k>k_0$. 
Hence $(v_k,u_k)$ converges in $C^{\infty}_{loc}(\Omega\times ]-T,0[\cup ]0,T[)$ to a smooth pair $(\bar{v},\bar{u})$ solving the equations \eqref{e:id-1} in $\R^n\times ]0,T[$ and such that \eqref{e:init(a)}, \eqref{e:init(b')}, \eqref{e:init(b)} and \eqref{e:init(c)} hold. It remains to show that $\frac{1}{2}|\bar{v}(x,0)|^2=1$ for almost every $x\in\Omega$.

>From \eqref{e:id-norm} we obtain
\begin{equation*}
\alpha_{k+1}\leq \alpha_k-\frac{1}{4}\min\bigl\{\alpha_k,C\alpha_k^2\bigr\}+|\Omega_{k+1}\setminus\Omega_k|\leq \alpha_k-\frac{1}{4}\min\bigl\{\alpha_k,C\alpha_k^2\bigr\} +2^{-k},
\end{equation*}
from which we deduce that
\begin{equation}\label{e:alfavan}
\alpha_k\to 0\textrm{ as }k\to\infty.
\end{equation}
Note that
\begin{equation}\label{e:id-norm2}
0\;\geq\;
\int_{\Omega}\Bigr[\frac{1}{2}|v_k(x,0)|^2-1\Bigl]\,dx
\geq -\bigl(\alpha_k+|\Omega\setminus\Omega_k|\bigr)
\;\geq\; -(\alpha_k + 2^{-k})\, .
\end{equation}
Therefore, by \eqref{e:alfavan},
\begin{equation}\label{e:intvan}
\lim_{k\uparrow\infty}
\int_{\Omega}\Bigr[\frac{1}{2}|v_k(x,0)|^2-1\Bigl]\,dx\;=\; 0\, .
\end{equation}
Finally, observe that, using \eqref{e:id-weak2}, for $t=0$ for every $k$
\begin{equation}\label{e:weakmolly}
\begin{split}
\Vert v_k*\rho_{\eta_k}-\bar{v}*\rho_{\eta_k}\Vert_{L^2}& \leq \sum_{j=0}^{\infty}\Vert v_{k+j}*\rho_{\eta_k}-v_{k+j+1}*\rho_{\eta_k}\Vert_{L^2}\\
&\leq 2^{-k}+2^{-(k+1)}+\dots\, \leq 2^{-(k-1)}
\end{split}
\end{equation}
and on the other hand
\begin{equation*}
\Vert v_k-\bar{v}\Vert_{L^2}\leq \Vert v_k-v_k*\rho_{\eta_k}\Vert_{L^2}+\Vert v_k*\rho_{\eta_k}-\bar{v}*\rho_{\eta_k}\Vert_{L^2} +\Vert \bar{v}*\rho_{\eta_k}-\bar{v}\Vert_{L^2}.
\end{equation*}
Thus, \eqref{e:strongmolly} and \eqref{e:weakmolly} imply that
$v_k(\cdot,0)\to \bar{v}(\cdot,0)$ strongly in $L^2(\R^n)$, which 
together with \eqref{e:intvan} implies that
$$
\frac{1}{2}|\bar{v}(x,0)|^2=1\textrm{ for almost every }x\in\Omega.
$$
\end{proof}

%%%%%%%%%%%%%%%%%%%%%%%%%%%%%%%%%%%%%%%%%

\section{Proofs of Theorem \ref{t:main} and Theorem \ref{t:psistema}}

Before embarking on the proof of Theorems \ref{t:main} and \ref{t:psistema}, we recall the following
well known fact: in the class of weak solutions in $C([0,T];L^2_w)$ it is possible to "glue" solutions which agree at a certain time. This is a consequence of the fact that if $v\in  C([0,T];L^2_w)$, then being a solution of \eqref{euler} in $\R^n\times [0,T[$ in the sense of distributions is {\em equivalent} to 
$v$ being divergence--free for all $t\in[0,T]$ and
\begin{equation*}
\int_{\R^n}v(x,t)\varphi(x,t)dx-\int_{\R^n}v(x,s)\varphi(x,s)dx=\int_s^t\int_{\R^n}[v\partial_t\varphi+\langle v\otimes v,\nabla\varphi\rangle]\,dx\,d\tau
\end{equation*}
for all $\varphi\in C_c^{\infty}(\R^n\times [0,T])$ with $\div\varphi=0$ and for all $s,t\in [0,T]$.

\subsection{Theorem \ref{t:main}}
{\bf Proof of (a)} Let $T=1/2$, $\Omega$ be the open
unit ball in $\R^n$, and $(\bar{v}, \bar{u})$ be as in
Proposition \ref{p:initialdata}. Define
$\bar{e}\equiv1$, $q_0\equiv0$,
\begin{equation}\label{e:defv0}
v_0 (x,t) \;:=\;
\left\{\begin{array}{ll}
\bar{v} (x,t)& \mbox{for $t\in [0,1/2]$}\\
\bar{v} (x, t-1)& \mbox{for $t\in [1/2,1]$,}
\end{array}\right.
\end{equation}
\begin{equation}\label{e:defu0}
u_0 (x,t) \;:=\;
\left\{\begin{array}{ll}
\bar{u} (x,t)& \mbox{for $t\in [0,1/2]$}\\
\bar{u} (x, t-1)& \mbox{for $t\in [1/2,1]$.}
\end{array}\right.
\end{equation}
It is easy to see that 
the triple $(v_0, u_0, q_0)$ 
satisfies the assumptions of Proposition \ref{p:main} with $\bar{e}\equiv 1$.
Therefore, there exists infinitely many solutions $v\in C ([0,1];L^2_w)$ 
of \eqref{euler} in $\R^n\times [0,1]$ with
$$
v (x, 0) = \bar{v}(x,0) = v(x, 1)
\textrm{ for a.e. }x\in\Omega,
$$
and such that
\begin{equation}\label{e:||=1}
\frac{1}{2}|v (\cdot ,t)|^2\;=\;{\mathbf 1}_\Omega \qquad
\mbox{for every $t\in ]0,1[$}\, .
\end{equation}
Since $\frac{1}{2}|v_0 (\cdot, 0)|^2={\mathbf 1}_\Omega$ 
as well, it turns out that the map $t\mapsto v (\cdot, t)$ is 
continuous in the strong topology of $L^2$.

By the remark above 
each such $v$ can be extended to a solution in $\R^n\times [0, \infty[$
which is 1-periodic in time, by setting $v (x, t)= v (x, t-k)$ for 
$t\in[k, k+1]$.
Then the energy 
$$
E (t) \;=\; \frac{1}{2} \int_{\R^n} |v(x,t)|^2 dx
$$
is equal to $|\Omega|$ at {\em every}
time $t$, and hence $v$ satisfies the strong energy equality
in the sense specified in Section 2. 

Next, notice that $\frac{1}{2}|v|^2 = {\mathbf 1}_{\Omega\times [0,\infty[}$ 
and that $p = - |v|^2/n = \textstyle{-\frac{2}{n}} \, {\mathbf 1}_{\Omega\times
[0, \infty[}$. Therefore for any $\varphi\in C_c^{\infty}(\R^n\times ]0,\infty[)$ we have
\begin{equation*}
\begin{split}
\int_{0}^{\infty}&\int_{\R^n}\frac{|v|^2}{2}\partial_t\varphi+  \left(\frac{|v|^2}{2}+p\right) v\cdot\nabla\varphi\,dx\,dt
\;=\\
& =\int_0^{\infty}\int_{\Omega} \partial_t\varphi +  \frac{n-2}{n} v\cdot \nabla\varphi\,dx\,dt\;=\;
\frac{n-2}{n}\int_0^{\infty}\int_{\R^n} v\cdot\nabla\varphi\,dx\,dt=0.
\end{split}
\end{equation*} 
This gives infinitely many solutions satisfying both the
strong energy equality and the local energy equality
and all taking the same initial data.

\medskip

{\bf Proof of (b)}  As in the proof of (a), let $T=1/2$, $\Omega$ be the open
unit ball in $\R^n$, and $(\bar{v}, \bar{u})$ be as in
Proposition \ref{p:initialdata}. Again, as in the proof of
(a) we set $q_0\equiv 0$. However we
choose $v_0$, $u_0$ and $\bar{e}$ differently: 
\begin{equation}\label{e:defv0b}
v_0 (x,t) \;:=\;
\left\{\begin{array}{ll}
\bar{v} (x,t)& \mbox{for $t\in [0,1/2]$}\\
0& \mbox{for $t\in [1/2,1]$,}
\end{array}\right.
\end{equation}
and
\begin{equation}\label{e:defu0b}
u_0 (x,t) \;:=\;
\left\{\begin{array}{ll}
\bar{u} (x,t)& \mbox{for $t\in [0,1/2]$}\\
0 & \mbox{for $t\in [1/2,1]$.}
\end{array}\right.
\end{equation}
Next consider the function 
$$
\tilde{e} (t)\;=\;\begin{cases}\max_{x\in \Omega} e (v_0 (x,t), u_0 (x,t))\,& 
\qquad \mbox{for $t\in ]0,1]$}\\
1&\qquad \mbox{for $t=0$}.\end{cases}
$$
It is easy to see that $\tilde{e}$ is continuous in $[0,1]$ (the continuity at $t=0$ follows 
from \eqref{e:init(d)} and \eqref{e:init(a)}), $\tilde e(t)<1$ for $t>0$, and $\tilde e=0$ in a neighborhood
of $t=1$. Define $\hat{e}: [0,1]\to \R$ as 
$$
\hat{e} (t) :=(1-t) + t \max_{\tau\in [t,1]} \tilde{e} (\tau).
$$
Then $\hat{e}$ is a continuous monotone decreasing
function, with 
$$
\hat{e} (0)=1,\,\hat{e} (1)=0,\textrm{ and }1>\hat{e} (t)>\tilde{e} (t)\textrm{ for every $t\in ]0,1[$}.
$$

Now, apply Proposition \ref{p:main} to get solutions
$v\in C ([0,1]; L^2_w)$ of \eqref{euler} in $\R^n\times [0,T]$ with
$v (\cdot, 0) = v_0 (\cdot, 0)$, $v (\cdot, 1) =0$ and such that
\begin{equation}\label{e:||=1c}
\frac{1}{2}|v (\cdot ,t)|^2\;=\;\hat{e} (t) \, {\mathbf 1}_\Omega \qquad
\mbox{for every $t\in ]0,1[$}\, .
\end{equation}
Arguing as in the proof of (a), we conclude that
$t\mapsto v (\cdot, t)$ is a strongly continuous map.
Since $v(\cdot ,1)=0$, we can extend $v$ by zero
on $\R^n\times [1,\infty[$ to get a global weak
solution on $\R^n\times [0,\infty[$. Clearly, this solution satisfies the strong
energy inequality. However, it does not satisfy the energy
equality. Note, in passing, that $v$ satisfies the local
energy inequality by the same reason as in (a).

\medskip

{\bf Proof of (c)} As in the proof of (a) and (b), 
let $T=1/2$, $\Omega$ be the open
unit ball in $\R^n$, and $(\bar{v}, \bar{u})$ be as in
Proposition \ref{p:initialdata}. Again, as in the proof of
(a) and (b) we set $q_0\equiv 0$. This time we
choose $v_0$, $u_0$ as in (b) and $\bar{e}$ as in (a).

Let $v_1\in C([0,1];L^2_w)$ be a solution of \eqref{euler}
obtained in Proposition \ref{p:main}. 
Since $\frac{1}{2}|v_0 (\cdot, 0)|^2={\mathbf 1}_\Omega$,
as before, the map $t\mapsto v_1(\cdot, t)$ is 
continuous in the strong topology of $L^2$ at every $t\in [0,1[$.
However, this map is {\em not} strongly continous at $t=1$,
because $v_1 (1, \cdot)=0$.  

Next, let $v_2\in C([0,1];L^2_w)$ be a solution of \eqref{euler} obtained
in Proposition \ref{p:main} with $\bar{e}\equiv 1$ and $(v_0,u_0,q_0)\equiv (0,0,0)$.
Since $v_1,v_2\in C([0,1];L^2_w)$ with $v_1(\cdot ,1)=v_2(\cdot ,0)=v_2(\cdot, 1)=0$,
the velocity field $v:\R^n\times [0, \infty[\,\to\R^n$ defined by
\begin{equation}\label{e:defc}
v (x,t) \;=\;
\left\{\begin{array}{ll}
v_1 (x, t)&\mbox{for $t\in [0,1]$}\\
v_2 (x, t-k)&\mbox{for $t\in [k, k+1], k=1,2,\dots$}
\end{array}\right.
\end{equation} 
belongs to the space $C([0,\infty[\,;L^2_w)$ and 
therefore $v$ solves \eqref{euler}. Moreover
$$
\frac{1}{2}\int |v(x,t)|^2\, dx \;=\; |\Omega|
\qquad\mbox{for every $t\not\in \N$} 
$$
and
$$
\frac{1}{2}\int |v(x,t)|^2\, dx \;=\; 0
\qquad\mbox{for every $t\in \N$, $t\geq 1$.} 
$$
Hence $v$ satisfies the weak energy inequality
but not the strong energy inequality.

\subsection{Theorem \ref{t:psistema}} We recall
that $p (\rho)$ is a function with $p' (\rho)>0$. 
Let 
$$
\alpha := p (1),\,\beta := p(2)\,\textrm{ and }\gamma = \beta-\alpha.
$$
Let $\Omega$ be the unit ball. Arguing
as in the proof of Theorem \ref{t:main}(a) we find
an initial data $v^0\in L^\infty(\R^n)$ 
with $|v^0|^2 = n \gamma \, {\mathbf 1}_\Omega$ and for
which there exist infinitely many weak solutions
$(v, \tilde{p})$ of \eqref{euler} with the
following properties:
\begin{itemize}
\item $v\in C ([0, \infty[; L^2)$ and
$|v|^2 = n\gamma \, {\mathbf 1}_{\Omega\times [0, \infty[}$;
\item $\tilde{p} = - |v|^2/n = - \gamma \, {\mathbf 1}_{\Omega\times [0, \infty[}$. 
\end{itemize}
In particular $v$ is divergence--free and $(v,\tilde p)$ satisfy
$$
\partial_tv+\div v\otimes v+\nabla \tilde{p}=0\quad\textrm{ in }\mathcal{D}'(\R^n\times ]0,\infty[).
$$
Then $(v,\hat p)$ also satisfy this equation, where $\hat p(x,t):=\tilde p(x,t)+\beta$. But observe that
for every $t\geq 0$ and for almost every $x\in\R^n$ we have
$$
\hat{p}(x,t)=\begin{cases}\alpha&\textrm{ if }x\in\Omega,\\
\beta&\textrm{ if }x\notin\Omega.\end{cases}
$$
so that
$$
\hat{p}(x,t)=p(\rho(x,t))\,\textrm{ for a.e. }(x,t)\in\R^n\times [0,\infty[,
$$
where $\rho$ is defined by
$$
\rho(x,t)=\begin{cases}1&\textrm{ if }x\in\Omega,\\
2& \textrm{ if }x\notin\Omega,\end{cases}
$$
for every $t\geq 0$. This shows that \eqref{e:test2} holds. To see that \eqref{e:test1} holds, observe that $\rho$ is independent of $t$ and $v$ is supported in $\Omega$. Hence
\begin{equation*}
\begin{split}
\int_0^{\infty}&\int_{\R^n}\bigl[\rho\partial_t\psi+\rho v\cdot\nabla\psi\bigr]\,dx\,dt+\int_{\R^n}\rho^0(x)\psi(x,0)\,dx=\\
&=\int_0^{\infty}\int_{\Omega}v\cdot\nabla\psi\,dx\,dt=\int_0^{\infty}\int_{\R^n}v\cdot
\nabla\psi\,dx\,dt=0,
\end{split}
\end{equation*}
because $v$ is divergence--free in $\R^n$ for every $t$.
Therefore for any such $v$, the pair $(\rho, v)$ 
is a weak solution of \eqref{e:psistema} with
initial data $(\rho^0, v^0)$, where $\rho^0={\mathbf 1}_\Omega+2\, {\mathbf 1}_{\R^n\setminus
\Omega}$. 

Each such solution is admissible. Indeed, similarly to the previous calculation we obtain
\begin{equation*}
\begin{split}
&\int_0^{\infty}\int_{\R^n}\bigl(\rho\eps(\rho)+\rho\frac{|v|^2}{2}\bigr)\partial_t\psi+\bigl(\rho\eps(\rho)+\rho\frac{|v|^2}{2}+p(\rho)\bigr)v\cdot\nabla\psi\,dx\,dt+\\
&\qquad+\int_{\R^n}\bigl(\rho^0\eps(\rho^0)+\rho^0\frac{|v^0|^2}{2}\bigr)\psi(x,0)\,dx=\\
&=\int_0^{\infty}\int_{\Omega}\frac{|v|^2}{2}\partial_t\psi+(\eps(1)+n\gamma+\alpha)v\cdot\nabla\psi
\,dx\,dt+\int_{\Omega}\frac{|v^0|^2}{2}\psi(x,0)\,dx\\
&=0,
\end{split}
\end{equation*}
because $v\in C([0,\infty[;L^2_w)$ and $v$ is divergence-free in $\R^n$.
This proves \eqref{e:entropy2} and thus concludes the proof of the theorem.

\section{Appendix A: Weak continuity in time for evolution equations}

In this section we prove a general lemma on the weak continuity
in time for certain evolution equations. Lemma \ref{l:weak}
is a corollary of this Lemma and standard estimates for the Euler and
Navier--Stokes equations.

\begin{lemma}\label{l:weakcont} Let $v\in L^\infty (]0,T[; L^2 (\R^n))$,
$u\in L^1_{loc} (\R^n\times ]0,T[, \R^{n\times n})$ and $q\in L^1_{loc}
(]0,T[\times \R^n)$ be distributional solutions of 
\begin{equation}\label{e:distrib}
\partial_t v + {\rm div}_x u + \nabla q\;=\; 0\, .
\end{equation}
Then, after redefining $v$ on a set of $t$'s of measure zero, $v\in C (]0,T[; L^2_w)$. 
\end{lemma}
\begin{proof} Consider a countable set $\{\varphi_i\}\subset C^\infty_c (\R^n, \R^n)$
dense in the strong topology of $L^2 $. Fix $\varphi_i$ and
any test function $\chi\in C^\infty_c (]0,T[)$. Testing \eqref{e:distrib}
with $\chi (t)\varphi_i (x)$ we obtain the following identity:
\begin{equation}\label{e:tested}
\int_0^T \Phi_i \partial_t \chi
\;=\; -\int_0^T \chi \int_{\rn{n}} \big[\langle u , \nabla \varphi_i\rangle +
q\, {\rm div} \varphi_i\big]\, ,
\end{equation}
where $\Phi_i (t):= \int \varphi_i (x) \cdot v (x,t) dx$. We conclude therefore that
$\Phi_i'\in L^1$ in the sense of distributions. Hence we can redefine each $\Phi_i$
on a set of times $\tau_i\subset ]0,T[$ of measure zero in such a way that $\Phi_i$ is continuous. We keep the same notation for these functions, and let $\tau=\cup_i\tau_i$. 
Then $\tau\subset ]0,T[$ is of measure zero and for every $t\in ]0,T[\setminus \tau$ we have
\begin{equation}\label{e:ae}
\Phi_i (t)\;=\; \int \varphi_i (x) \cdot v (x,t)\, dx
\qquad \mbox{for every $i$.}
\end{equation}
Moreover, with $c:= \|v\|_{L^\infty_t (L^2_x)}$ we have that
$|\Phi_i (t)|\leq c \|\varphi_i\|_{L^2}$ for all $t\in]0,T[$. Therefore, for each $t\in]0,T[$
there exists a unique
bounded linear functional $L_t$ on $L^2 (\rn{n}, \rn{n})$ such that $L_t (\varphi_i)
= \Phi_i (t)$. By the Riesz representation theorem there exists 
$\bar{v}(\cdot,t)\in L^2(\R^n)$ such that 
\begin{itemize}
\item $\bar{v} (\cdot ,t)=v(\cdot,t)$ for every $t\in ]0,T[\setminus\tau$;
\item $\|\bar{v} (\cdot ,t)\|_{L^2}\leq c$ for every $t$;
\item $\int \bar{v} (x,t) \cdot \varphi_i (x) dx = \Phi_i (t)$ for every $t$.
\end{itemize}
To conclude we show that $\bar{v}\in C(]0,T[;L^2_w)$, i.e.~that
for any $\varphi\in L^2(\R^n,\R^n)$ the function
$\Phi (t):= \int v (x,t)\cdot \varphi (x) dx$ is continuous on $]0,T[$.
Since the set $\{\varphi_i\}$ is dense in $L^2(\R^n,\R^n)$, we can find
a sequence sequence $\{j_k\}$ such that $\varphi_{j_k}\to \varphi$ strongly in $L^2$.
Then
\begin{equation}\label{e:est}
|\Phi (t) - \Phi_{j_k} (t)| \;\leq\; c \|\varphi_{j_k} - \varphi\|_{L^2}\, .
\end{equation}
Therefore $\Phi_{j_k}$ converges uniformly to $\Phi$, from which we
derive the continuity of $\Phi$. This shows that $\bar{v}\in C(]0,T[;L^2_w)$ and
concludes the proof.
\end{proof}

\section{Appendix B: Dissipative solutions}

We follow here the book \cite{Lions} and define dissipative solutions
of \eqref{euler}. First of all, for any divergence--free vector
field $v\in L^2_{loc} (\R^n \times [0,T])$ we consider the following
two distributions:
\begin{itemize}
\item The symmetric part of the gradient $d(v):=\frac{1}{2} (\nabla v + \nabla v^t)$;
\item $E(v)$ given by 
\begin{equation}\label{e:E}
E(v)\;:=\; -\partial_t v - P ({\rm div}\, (v\otimes v))\, .
\end{equation}
\end{itemize}
Here $P$ denotes the Helmholtz projection on divergence--free fields, so that
if $p(x,t)$ is the potential--theoretic solution of $- \Delta p 
= \sum_{i,j} \partial^2_{ij} (v^iv^j)$, then 
$$
P ({\rm div}\, (v\otimes v)) 
\;=\; {\rm div}\, (v\otimes v)  + \nabla p \, .
$$
Finally, when $d(v)$ is locally summable, we denote by $d^- (v)$    
the negative part of its smallest eigenvalue, 
that is $(-\lambda_{min} (d(v)))^+$. 

P.~L.~Lions introduced the following definition in \cite{Lions}:

\begin{definition}\label{d:dissipative}
Let $v\in L^\infty ([0,T]; L^2 (\R^n))\cap C([0,T]; L^2_w)$. Then $v$
is a dissipative solution of \eqref{euler} if the following two conditions hold
\begin{itemize}
\item $v(x, 0)=v_0(x)$ for $x\in \R^n$; 
\item ${\rm div}\, v = 0$ in the sense of distributions;
\item whenever $w\in C ([0,T]; L^2 (\R^n))$ is such that
$d(w)\in L^1_t(L^\infty_x)$, $E(w)\in L^1_t(L^2_x)$ and $\textrm{div }w=0$, then
\begin{equation}\label{e:dissipineq}
\begin{split}
\|v& (\cdot ,t)-w (\cdot ,t)\|^2_{L^2_x}\leq
e^{\int_0^t 2\|d^- (w)\|_{L^\infty_x} d\tau}  \|v_0(\cdot) - w (\cdot, 0)\|^2_{L^2_x}\\
&+ 2\int_0^t\int_{\R^n}  e^{\int_s^t 2\|d^- (w) \|_{L^\infty_x} d\tau} E(w) (x,s)
\cdot(v(x,s)-w(x,s))\, dx\, ds
\end{split}
\end{equation}
for every $t\in [0,T]$. 
\end{itemize}
\end{definition}

We next come to the proof of Proposition \ref{p:WImpDissBis}
which we state again for the reader's convenience.

\begin{proposition}\label{p:WImpDiss}
Let $v\in C ([0, T]; L^2_w)$ be a weak solution of \eqref{euler}
satisfying the weak energy inequality. Then $v$ is a dissipative solution.
\end{proposition}

\begin{proof} As already remarked at page 156 of \cite{Lions} it suffices
to check Definition \ref{d:dissipative} for smooth $w$. This is achieved by
suitably regularizing the test function $w$ of \eqref{e:dissipineq} and observing
that if $w\in C ([0,T]; L^2 (\R^n))$ is such that
$d(w)\in L^1_t(L^\infty_x)$, then any approximation $w_k$ such that 
\begin{enumerate}
\item[(a)] $w_k\to w$ in $C([0,T];L^2)$;
\item[(b)] $d(w_k)\to d(w)$ a.e. in $\R^n\times [0,T]$;
\item[(c)] $\limsup_{k\to\infty}\|d (w_k)\|_{L^\infty_x}\leq \|d (w)\|_{L^\infty_x}$
\end{enumerate}
also satisfies
$$
E(w_k)\to E(w)\textrm{ in }L^1_tL^2_x
$$
and hence one can pass to the limit in \eqref{e:dissipineq}. Indeed, this follows from the observation that $P(E(w))= 2 P(d (w)\cdot w)$ (see the computations on page 155 of \cite{Lions}).

\medskip

{\bf Step 1.} Next we show that it suffices to check 
Definition \ref{d:dissipative}
when $w$ is compactly supported in space. Indeed, fix $w$ 
as above. We claim that we can approximate
$w$ with compactly supported divergence--free vector fields 
$w_k$ such that (a),(b) and (c) above hold. The reader may 
consult Appendix A of \cite{Lions} and jump directly to
Step 2. Otherwise, the following is a short self-contained
proof.

Fix a smooth cut--off function
$\chi$ equal to $1$ on the ball $B_1 (0)$, supported in the ball $B_2 (0)$,
and taking values between $0$ and $1$, and set $\chi_r (x)= \chi (r^{-1} x)$. 
Let $\xi$ be the potential--theoretic solution of $\Delta \xi = \textrm{curl }w$, so that $w=\textrm{curl }\xi$. Recall that in dimension $n=2$ the curl operator can be defined as $\textrm{curl }=(-\partial_2,\partial_1)$, in dimension $n=3$ it is given by $\textrm{curl }w=\nabla\times w$ and $\xi$ is obtained via the Biot-Savart law. Let $\langle\xi\rangle_k=\frac{1}{|B_{2k}\setminus B_k|}\int_{B_{2k}\setminus B_k}\xi\,dx$ and let
$$
w_k = \textrm{curl} \bigl(\chi_{k} (\xi-\langle\xi\rangle_k)\bigr).
$$
Clearly $w_k$ is compactly supported and divergence--free. Since $\xi$ is smooth, and $\|\partial_i (\chi_{k})\|_{\infty}\leq C k^{-1}$ and $\|\partial^2_{ij} 
(\chi_{k})\|_{\infty}\leq C k^{-2}$, we see that 
$$
d(w_k)(\cdot,t)\to d(w)(\cdot,t)\textrm{ locally uniformly }
$$
for every $t$. Thus (b),(c) follow easily.
Moreover $\|\nabla \xi(\cdot,t)\|_{L^2_x}\leq \|w(\cdot,t)\|_{L^2_x}$ and hence, using the Poincar\'e inequality, for every $t\in[0,T]$ we have
\begin{equation*}
\begin{split}
\|w_k-w\|^2_{L^2_x}&\leq C\int_{\R^n\setminus B_k (0)} |w|^2dx
+ C \|\nabla \chi_{k^{-1}}\|_{C^0}^2 \int_{B_{2k} (0)\setminus B_k (0)}
|\xi-\langle\xi\rangle_k|^2dx\\
&\leq C\int_{\R^n\setminus B_k (0)} |w|^2dx +
\frac{C}{k^2} \int_{B_{2k}\setminus B_k (0)} |\nabla\xi|^{2}dx\\
&\leq C\int_{\R^n\setminus B_k (0)} |w|^2dx +
\frac{C}{k^2} \int_{\R^n} |w|^{2}dx.
\end{split}
\end{equation*}
Since $w\in C([0,T];L^2(\R^n))$, we deduce (a).

\medskip

{\bf Step 2.} We are now left with task of showing \eqref{e:dissipineq}
when $w$ is a smooth test function compactly supported in space. 
Consider the function
$$
F(t) := \int_{\R^n} |w(x,t)-v(x,t)|^2 dx.
$$
Since $w$ is smooth
and $v\in C ([0,T], L^2_w)$, $F$ is 
lower--semicontinuous. Moreover, due to the weak
energy inequality $v(t, \cdot) \to v(0, \cdot)$ strongly
in $L^2_{loc}$ as $t\downarrow 0$. So $F$ is continuous
at $0$. We claim that, in the sense of distributions,
\begin{equation}\label{e:IneqForm}
\frac{dF}{dt}\leq
2 \int_{\R^n} \bigl[E(w) \cdot (v-w)-d(w) (v-w)\cdot (v-w)\bigr]\, dx\, .
\end{equation}
>From this inequality we infer 
\begin{equation}\label{e:ineq2}
\frac{dF}{dt} \;\leq\; 2 \|d^- (w) (t, \cdot)\|_{L^\infty_x} F(t)
+  2 \int_{\R^n} \big[E(w) \cdot (v-w)\big]\, dx\, .
\end{equation}
>From the continuity of $F$ at $t=0$ and Gronwall's Lemma,
we conclude \eqref{e:dissipineq} for a.e. $t$. By the lower
semicontinuity of $F$, \eqref{e:dissipineq} actually
holds for every $t$. Therefore it remains to prove \eqref{e:IneqForm}.
We expand $F$ as
\begin{equation*}
\begin{split}
F(t) &= \int_{\R^n} |v(x,t)|^2 \, dx
+ \int_{\R^n} |w(x,t)|^2 \, dx - 2 \int_{\R^n} \big[v(x,t)\cdot w(x,t)\big]\,dx\\
&=: F_1 (t) + F_2 (t) + F_3 (t)\, .
\end{split}
\end{equation*}
The weak energy inequality implies $\frac{d}{dt} 
F_1 (t)\leq 0$ and a standard calculation gives
$$
\frac{d F_2}{dt} (t)\;=\; - 2 \int_{\R^n} \big[ E(w)\cdot w\big]\, dx\, .
$$
It remains to show that
\begin{equation}\label{e:IneqForm2}
\frac{dF_3}{dt}=
2 \int_{\R^n} \bigl[E(w) \cdot v-d(w)(v-w)\cdot (v-w)\bigr]\, dx
\end{equation}
We fix a smooth function $\psi\in C^\infty_c (]0,T[)$ 
and test \eqref{euler} (or more precisely \eqref{distrib}) with $w(x,t)\psi (t)$. It then 
follows that
\begin{equation}\label{e:ineq10}
2 \int_{\R}\int_{\R^n} v\cdot w \psi'\, dx\, dt
=- 2 \int_{\R}\psi\int_{\R^n} \bigl[v\cdot \partial_t w + 
\langle v\otimes v , \nabla w\rangle\bigr]\,dx\,dt .
\end{equation}
Inserting $\partial_t w = - E(w) - P ({\rm div} 
(w\otimes w))$ and taking into account that ${\rm div}\, v
=0$, we obtain
\begin{eqnarray}
\int_{\rn{}} F_3 (t) \psi' (t)\, dt
&=& 2 \int_{\rn{}} \psi \int_{\rn{n}}
\bigl[ \langle v\otimes v , \nabla w\rangle
- {\rm div} (w\otimes w)\cdot v\bigr]\, dx\,dt\nonumber\\
&&- 2\int_{\rn{}} \psi \int_{\rn{n}} E(w) \cdot v\, dx\,dt
\label{e:ineq11}
\end{eqnarray}
Next, observe that ${\rm div} (w\otimes w)\cdot v
= \sum_{j,i} v_j w_i \partial_i w_j$ and
that $\langle v\otimes v, \nabla w\rangle
= \sum_{j,i} v_j v_i \partial_i w_j$. Therefore we have
\begin{equation}\label{e:eq12}
\langle v\otimes v, \nabla w\rangle
- {\rm div} (w\otimes w)\cdot v\;=\;
\nabla w \,(v-w)\cdot v  \, .
\end{equation}
On the other hand, 
$$
\nabla w\, (v-w)\cdot w
\;=\; \sum_{i,j} (v_i-w_i) \partial_i w_j w_j
\;=\; (v-w)\cdot \nabla \frac{1}{2}|w|^2\, .
$$
Since $v-w$ is divergence--free in the sense
of distributions and $|w|^2/2$ is a smooth function
compactly supported in space,  
integrating by parts we get
\begin{equation}\label{e:eq13}
\int_{\rn{}} \psi \int_{\rn{n}}
\bigl[\nabla w\,(v-w)\cdot w\bigr] \, dx\,dt \;=\; 0\, .
\end{equation}
>From \eqref{e:ineq11}, \eqref{e:eq12} and \eqref{e:eq13}
we obtain
\begin{eqnarray}
\int_{\rn{}} F_3 (t) \psi' (t)\, dt
&=& 2\int_{\rn{}} \psi \int_{\rn{n}}
\bigl[\nabla w\,(v-w)\cdot (v-w)\bigr]
\, dx\,dt\nonumber\\
&&- 2 \int_{\rn{}} \psi \int_{\rn{n}}
E(w) \cdot v\, dx\,dt\, .
\label{e:ineq14}
\end{eqnarray}
Finally, observe that 
$$
\nabla w (v-w)\cdot (v-w)
= \bigl\langle \nabla w , (v-w)\otimes (v-w)\bigr\rangle=
\bigl\langle d(w) , (v-w)\otimes (v-w)\bigr\rangle,
$$
since $(v-w)\otimes (v-w)$ is a symmetric matrix. 
Plugging this into \eqref{e:ineq14}, by the
arbitrariness of the test function $\psi$, we obtain
\eqref{e:IneqForm2}.
\end{proof}

\bibliographystyle{acm}

\end{document}